\documentclass{amsart}
\usepackage[USenglish]{babel}
\usepackage{color}
\usepackage{graphicx}

\usepackage{amssymb}
\usepackage{amsmath}

\newcommand{\R}{\mathrm{I\!R\!}}

\newtheorem{theorem}{Theorem}[section]
\newtheorem{lemma}[theorem]{Lemma}
\newtheorem{proposition}[theorem]{Proposition}

\newtheorem{remark}[theorem]{Remark}

\numberwithin{equation}{section}

\title[A Model for Tumor Invasion ]{An Analysis of a Mathematical Model Describing Acid-mediated Tumor Invasion}

\author{Anderson L. A. de Araujo}
\address{Departamento de Matem\'atica, Universidade Federal de Vi\c{c}osa, Vi\c{c}osa, MG, Brazil}
\email{anderson.araujo@ufv.br}
\thanks{The first author was was partially sponsored by FORTIS POS-GRADUA\c{C}\~AO EM MATEM\'ATICA - 9389}

\author{Artur C. Fassoni}
\address{Instituto de Matem\'atica e Computa\c{c}\~ao, Universidade Federal de Itajub\'a, Itajub\'a, MG, Brazil}
\email{fassoni@unifei.edu.br}

\author{Lu\'is F. Salvino }
\address{Instituto de Ci\^encias Exatas, Universidade Federal de Vi\c{c}osa - Campus Florestal, Florestal, MG, Brazil}
\email{luisfsalvino@gmail.com}



\begin{document}

\begin{abstract}
We present a mathematical analysis of a reaction-diffusion model describing acid-mediated tumor invasion. The model describes the spatial distribution and temporal evolution of tumor cells, normal cells, and excess lactic acid concentration. The model assumes that tumor-induced alteration of microenvironmental pH provides a simple but complete mechanism for cancer invasion. We provide results on the existence and uniqueness of a solution considering Neumann and Dirichlet boundary conditions. We also provide numerical simulations to the solutions considering both boundary conditions.


\end{abstract}

\maketitle


\par
\section{Introduction}

The main objective of this work is to perform a rigorous mathematical analysis of a system of nonlinear partial differential equations corresponding to a generalization of a mathematical model describing the growth of a tumor proposed in \cite{Fassoni}.

To describe the model,
let $\Omega \subset \R^{2}$, be an open and bounded set; let also  $0< T< \infty$ be a given final time  of interest and denote $t$ the times between $[0,T]$ and
$Q=\Omega\times(0,T)$, the space-time cylinder and
$\bar{\Gamma}=\partial \Omega\times(0,T)$,  the space-time boundary. Then, the system of equations we are considering is the following:
\begin{equation}
\label{0riginalEquations}
\left\{
\begin{array}{lcl}
\displaystyle
\frac{\partial N}{\partial t} = r_N - \mu_N N - \beta_1 N A - \alpha_H\gamma_H N H, & \textup{in}& Q,
\vspace{0.2cm}
\\
\displaystyle
\frac{\partial A}{\partial t} = \xi_A \Delta A + r_A A\left(1-\frac{A}{k_A}\right)-(\mu_A+\epsilon_A)A - \beta_3 NA , &\textup{in}& Q,
\vspace{0.2cm}
\\
\displaystyle
\frac{\partial H}{\partial t} = \xi_H \Delta H + \nu A - \tau_H H - \gamma_H NH , &\textup{in}& Q,
\vspace{0.2cm}
\\
\displaystyle \frac{\partial A}{\partial \eta} (\cdot) =\frac{\partial H}{\partial \eta} (\cdot) =0, &\textup{on}& \Gamma,
\vspace{0.2cm}
\\
\displaystyle N(\cdot,0) = N_0(\cdot), A(\cdot,0) = A_{0}(\cdot), H(\cdot,0) = H_0(\cdot), &\textup{in}& \Omega.
\end{array}
\right.
\end{equation}
In \cite{Fassoni}, Fassoni studied a simplified model given by a system of ODE's describing the growth of a tumor and its effect in the normal tissue, together with the tissue response to the tumor and the effect of chemotherapeutic treatments. The aim of the authors was to provide some $insights$ on the description of cancer onset and treatment as transitions between alternative basins of attractions. The model studied in \cite{Fassoni} is given by the following ODE system:
\begin{equation}
\label{início1}
\left\{
\begin{array}{lcl}
\displaystyle	
\frac{d N}{d t} = r_N - \mu_N N - \beta_1 N A - \alpha_N\gamma_N D N,
\vspace{0.2cm}
\\
\displaystyle	
\frac{d A}{d t} = r_A A\left(1-\frac{A}{k_A}\right)-(\mu_A+\epsilon_A)A - \beta_3 NA - \alpha_A\gamma_A D A,
\vspace{0.2cm}
\\
\displaystyle	
\frac{d D}{d t} = \mu - \gamma_A D A - \gamma_N D N - \tau D,
\end{array}
\right.
\end{equation}
where $N$ represents the normal cells in a given tissue of the human body, $A$ represents the tumor cells in this tissue and $D$ represents the concentration of a drug used to treat such a tumor.

In the model, parameter $r_N$ represents the total constant reproduction of normal cells and $\mu_N$ their natural mortality. A constant flow for normal cells is considered in the vital dynamics and does not depend on the density, as is the case of logistical growth generally assumed, see \cite{Jb}. The reason is that in a normal and already formed tissue, the imperative dynamics is not the intraspecific competition of cells for nutrients, but the maintenance of a homeostatic state, through the natural replenishment of old cells, see \cite{Simons}.

In contrast, cancer cells have independence on growth signals and maintain their own growth program, as for example, a growing embryonic tissue, see \cite{Fedi}. Thus, a density-dependent growth is considered for tumor cells. Several growth laws could be used, such as Gompertz, generalized logistic, Von Bertanlanfy and others, see \cite{Sarapata}, but the logistic growth is chosen due to its simplicity. A natural mortality rate $\mu_A$ and an extra mortality rate $\epsilon_A$, due to apoptosis, see \cite{Danial}, are also included.

Parameters $\beta_1$ and $\beta_3$ encompass in a simplified way the many negative interactions between tumor cells and normal. Parameter $\beta_3$ encompasses the effects imposed on cancer cells by normal cells, such as competition for nutrients and the release of anti-growth and death signals. In the same way, parameter $\beta_1$ encompasses mechanisms developed by tumor cells that damage the normal tissue, such as competition, release of death signals and increased local acidity.

The hypotheses that lead to the third equation of (\ref{início1}), describing the pharmacokinetics of the chemotherapeutic drug are the following. The drug has a clearance rate $\tau$. The rates of absorption and deactivation of the drug by normal and cancer cells are described in terms of the law of mass action with rates $\gamma_N$ and $\gamma_A$. Following the linear hypothesis of \cite{Andre}, the amounts of drug absorbed by normal ($\gamma_N N D$) and cancer cells ($\gamma_A A D$) kill such cells at rates $\alpha_N$ and $\alpha_A$, respectively. Finally, parameter $\mu$ represents a constant infusion rate mimicking a metronomic dosage, i.e., a near continuous and long-term administration of the drug.

System (\ref{início1}) is similar to the classic Lotka-Volterra competition model, commonly used in models for tumor growth and biological invasions, but there is a fundamental difference: the use of a constant flux instead of a logistic growth for normal cells breaks the symmetry observed in the classic Lotka-Volterra model, so that there will be no equilibrium for $N = 0$. Thus, normal cells will never be extinct, unlike these models. In \cite{Fassoni} the authors believe that this is not a problem, but on the contrary, it is a realistic result. In fact, roughly speaking, cancer does "don not win" by the fact that it kills all the cells in the tissue, but because it reaches a dangerous size that disrupts the proper functioning of the tissue and threatens the health of the individual. A constant flow term has already been taken in other well-known models of cancer, specifically, to describe the growth of immune cells, see \cite{Earn}.

In this work, we are not interested in analyzing the dynamics of model (\ref{início1}), since in \cite{Fassoni} a study has already been made on the equilibrium points and questions about stability or instability of such points. Our objective is to study the existence and uniqueness of solutions to the modified system (\ref{0riginalEquations}).

System (\ref{0riginalEquations}) considers the dynamics of normal and tumor cells under the same hypothesis as the system (\ref{início1}), and explicitly considers the production of lactic acid by cancer cells. Such acid production is well-known as a mediator of tissue invasion by tumor cells, and is a by-product of the altered metabolism of tumor cells, which exert glucose metabolism not by oxidative phosphorylation (as normal cells), but mainly by glycolysis. This altered metabolism is known as Warburg effect, see \cite{Jb}. The increased acid concentration in the tissue damages the normal cells, thereby destroying the ``barriers'' for tumor growth and contributing to tumor persistence and further tissue invasion.

Model (\ref{0riginalEquations}) describes such phenomena by considering the original equations for the normal cells ($N(t)$) and cancer cells ($A(x,t)$), and an additional equation for the concentration of extracellular lactic acid in excess of normal tissue acid concentrations ($H(x,t)$), produced at a rate $\nu$ by tumor cells, cleared by tissue vasculature at a rate $\tau_H$ and absorbed by normal cells at a rate $\gamma_H$, following the mass-action law. The absorbed amount $\gamma_H NH$ causes a damage to normal cells with rate $\alpha_H$. The aspect of acid-mediated tumor invasion clearly needs to be considered in a spatial setting, thereby it is better represented by a PDE model instead of an ODE system, which describes a spatially homogeneous tumor. Therefore, we assume a diffusive behavior for both the tumor cells and the extracellular lactic acid excess, with diffusion coefficients $\xi_A$ and $\xi_H$ respectively, while normal cells do not move. This is in agreement with previous models, see \cite{Jb}, and regards the additional capability of tumor cells to move within a ``rigid'' tissue. Previous works approached similar models, but all considering a density dependent growth for normal cells (logistic growth), see (\cite{Jb}) and references therein.

In this work, we are not interested in the treatment phase, but only on the invasion phase and the establishment of a tumor. In a future work, we may extend model (\ref{0riginalEquations}) by including a differential equation for the chemotherapeutic drug as given by system (\ref{início1}).

This paper is organized as follows. In Section 2 we present the technical hypothesis state our main result. In Section 3 we study an auxiliary problem. Using its solution, we prove our main result in Section 4. In Section 5 we present numerical simulations illustrating model behavior.

\section{Technical hypotheses and main result}
Recall $\Omega \subset \R^2$ be a domain with boundary $\partial\Omega$, $0\leq T <\infty$, and denote $Q = \Omega\times (0, T)$ and $\Gamma =  \partial \Omega\times (0, T)$.
We will use standard notations for Sobolev spaces, i.e., given $1~\leq~p~\leq~+\infty $ and $k \in\mathbb{N}$, we denote
$$W_{p}^{k}(\Omega)=\left\{ f \\ \in L^{p}(\Omega) : D^{\alpha}f \in L^{p}(\Omega), |\alpha| \leq k \right\}; $$
when $p=2$, as usual we denote $W_{2}^{k}(\Omega) = H^k (\Omega)$;
properties of these spaces can be found for instance in~Adams~\cite{Adams}.
Problem~(\ref{0riginalEquations}) will be studied in the standard functional spaces denoted by
\begin{eqnarray*}
  W_{q}^{2,1}(Q) & =& \left\{f\in L^{q}(Q):D^{\alpha}f\in L^{q}(Q), \, \forall
1\leq|\alpha|\leq 2,  f_t \in L^{q}(Q)\right\},
\end{eqnarray*}
\begin{eqnarray*}
  W &=& \left\{f\in L^\infty(Q): f_t \in L^\infty(Q)\right\}
\end{eqnarray*}
and
\begin{eqnarray*}
  L^{p}(0,T;B) &=&  \left\{f:(0,T)\rightarrow B:   \|f(t)\|_{L^{p}(0,T;B)} <+\infty \right\},
\end{eqnarray*}
where $B$ is suitable Banach space, and the norm is given by
$\|f(t)\|_{L^{p}(0,T;B)} = \|\ \|f(t)\|_{B}\ \|_{L^{p}((0,T))}$.
We remark that $L^{p}(Q) = L^{p}((0,T);L^{p}(\Omega))$.
Results concerning these spaces can be found  for instance in Ladyzhenskaya~\cite{Ladyzhenskaya}
and Mikhaylov~\cite{Mikhaylov}.

\vspace{0.1cm}
Next, we state some hypotheses that will be assumed throughout this article.

\subsection{Technical Hypotheses:}
\label{MainHypotheses}

\begin{itemize}

\item[{\bf (i)}] $\Omega\subset\mathbb{R}^2$ is a bounded $C^2$-domain.

\item[{\bf (ii)}] $0< T < \infty$, and $Q=\Omega\times(0,T)$.

\item[{\bf (iii)}] $N_0 \in L^{\infty}(\Omega)$ and $A_0, H_0 \in W^{\frac{3}{2}}_{4}(\Omega)$, satisfying $\frac{\partial A_0}{\partial \eta} (\cdot) =\frac{\partial H_0}{\partial \eta} (\cdot) =0, \textup{ on } \partial\Omega$.

\item[{\bf (iv)}] $0 \leq A_0 \leq k_A$ and $N_0, H_0 \ge 0$ a.e. on $\Omega$

\end{itemize}

\begin{remark}
The constraints imposed in~{\bf (iv)} on the initial conditions  are natural biological requirements.
\end{remark}

\subsection{Main result:}
\begin{theorem}
\label{Teorema1}
Assume that the Technical Hypotheses \ref{MainHypotheses} hold;
then, there exists a unique nonnegative solution $(N,A,H) \in W \times W^{2,1}_4(Q) \times W^{2,1}_4(Q)$ of Problem (\ref{0riginalEquations}). Moreover, $N, A$ and $H$ are functions satisfying
\begin{eqnarray*}
  N \leq ||N_0||_{L^\infty(Q)} + r_N T, \ A \leq k_A \ a.e. \ in \ Q
\end{eqnarray*}
and
\begin{eqnarray*}
  ||N||_W + ||A||_{ W^{2,1}_4(Q)} + ||H||_{ W^{2,1}_4(Q)} \leq C,
\end{eqnarray*}
where $C$ is a constant depending on $r_N$, $\mu_N$, $\beta_1$, $\alpha_H$, $\gamma_H$, $k_A$, $\nu$, $T$, $\Omega$, $||N_0||_{L^\infty (\Omega)}$, $||A_0||_{W^{\frac{3}{2}}_4(\Omega)}$, $||H_0||_{W^{\frac{3}{2}}_4(\Omega)}$ and $||H||_{L^\infty(Q)}$.
\end{theorem}

\begin{remark}
The explicit knowledge on how the constant $C$ appearing in the above estimates depends on the given data is important to application to related control problems.
\end{remark}

\subsection{Known technical results:}

To ease the references, we also state some technical results to be used in this paper.
The first one is consequence of Theorem~5.4, p. 97, in Adams~\cite{Adams}:
\begin{lemma}
\label{imersao_Sobolev}
Suppose that $\Omega\subset\mathbb{R}^{n}$ satisfies the cone property and $1 \leq p <\infty$.
The following continuous embeddings hold:
\[
\begin{array}{l}
\mbox{\bf (i)} \; \, W^{2}_{p}(\Omega)\hookrightarrow W^{2-\frac{2}{q}}_{q}(\Omega)$,
for all $p\leq q\leq 5p/3$, if $ n\leq 3;
\\
\mbox{\bf (ii)} $ If $ kp<n$, then $
\, W^{k}_{p}(\Omega)\hookrightarrow L^{q}(\Omega)$,
for all $p \leq q \leq np/(n-kp);
\\
\mbox{\bf (iii)} $ If $ kp=n$, then $
\, W^{k}_{p}(\Omega)\hookrightarrow L^{q}(\Omega)$,
for all $p \leq q <\infty;
\\
\mbox{\bf (iv)} $ If $ kp>n$, then $
\, W^{k}_{p}(\Omega)\hookrightarrow L^{\infty}(\Omega).
\end{array}
\]
\end{lemma}

The next result sometimes is called the Lions-Peetre embedding theorem (see Lions~\cite{Lions}, pp.15);
it is also a particular case of Lemma~3.3, pp.80, in Ladyzhenskaya~\cite{Ladyzhenskaya}:
(obtained by taking $l = 1$ and $r = s = 0$).
\begin{lemma}
\label{icontLp01}
Let $\Omega$ be a domain of $\R^n$ with boundary $\partial \Omega$ satisfying the cone property.
Then, the functional space $W^{2,1}_p(Q)$ is continuously embedded in $u \in L^{q}(Q)$ for $q$ satisfying: {\bf (i)} $1 \leq q \leq \frac{p(n+2)}{n+2-2p}$, if $ p< \frac{n+2}{2}$; {\bf (ii)} $1 \leq q <\infty$, if $p= \frac{n+2}{2}$ and {\bf (iii)} $q=\infty$, if $p>\frac{n+2}{2}$.
\noindent
In particular, for such $q$ and any function $u \in W^{2,1}_p(Q)$ we have that
\begin{eqnarray*}
\label{i.01}
\displaystyle \|u\|_{L^{q}(Q)} \leq C\|u\|_{W^{2,1}_p(Q)},
\end{eqnarray*}
\noindent
with a constant $C$ depending only on $\Omega$, $T$, $p$, $q$, $n$.

In the cases {\bf (ii)}, {\bf (iii)} or in {\bf (i)} when $\displaystyle 1 \leq q < \frac{p(n+2)}{n+2-2p}$, the referred embedding is compact.

\end{lemma}

\vspace{0.1cm}
Next, we consider the following general and simple parabolic initial-boundary value problem:
\begin{equation}
\label{P_Newmman}
\left\{
\begin{array}{lcl}
\displaystyle	\frac{\partial u}{\partial t} - \sum\limits_{i,j=1}^na_{ij}(x,t)\frac{\partial u^2}{\partial x_ix_j}
+ \sum\limits_{j=1}^na_i(x,t)\frac{\partial u}{\partial x_j} + a(x,t)u=f & \textup{in} & Q,
\vspace{0.2cm}
\\
\displaystyle
\sum\limits_{i=1}^n b_i(x,t)\frac{\partial u}{\partial x_i} + b(x,t)u =0
& \textup{on} & \Gamma,
\vspace{0.2cm}
\\
\displaystyle
u(\cdot,0)= u_0(\cdot) & \textup{in} & \Omega .
 \end{array}
 \right.
\end{equation}

Existence and uniqueness of solutions for this problem is a particular case of Theorem~$9.1$, pp.$341$,
in Ladyzenskaya~\cite{Ladyzhenskaya} for the case of Neumann boundary condition, according to the remarks at the end Chapter IV, section 9, p. 351 in \cite{Ladyzhenskaya}.
In the following, we state this particular result, stressing the dependencies certain norms of the coefficients,  that will be important in our future arguments.
\begin{proposition}
\label{sol. Neumann}
Let $\Omega$ be a bounded domain in $\mathbb{R}^n$, with a $C^{2}$ boundary $\partial \Omega$,
$a_{ij}$ be bounded continuous functions in $Q$,
and $q > 1$. Assume that
\begin{enumerate}

\item $a_{ij} \in C(\bar{Q})$, $i, j=1, \ldots, n$; $[a_{ij}]_{n \times n}$ is a real positive matrix such that for some positive constant $\beta$
we have  $ \sum\limits_{i,j=1}^n a_{ij}(x,t)\xi_i\xi_j\geq \beta|\xi|^2$ for all
$(x,t) \in Q$ and all $\xi \in R^n$,
;

\item $\displaystyle f \in L^p(Q)$;

\item $\displaystyle a_i \in L^r(Q)$ with
either $r =  \max\big(p, n + 2\big)$ if $p \neq n + 2$
or $r =  n + 2 + \varepsilon$, for any  $\varepsilon>0$, if $p = n + 2$;

\item $\displaystyle a \in L^s(Q)$ with
either $s = \max\big(p, (n + 2)/2\big)$ if  $\displaystyle p \neq (n + 2)/2$
or $s = (n + 2)/2 + \varepsilon$, for any $\varepsilon>0$, if $\displaystyle p = (n + 2)/2$.

\item $b_i, b \in C^2 (\bar{\Gamma})$, $i=1, \ldots, n$, and
the coefficients $b_i(x,t)$ satisfy the condition
$\left| \sum\limits_{i=1}^n b_i(x,t)\eta_i(x) \right|\geq \delta >0$
for $a.e.$ in  $\partial\Omega \times (0,T)$,
where $\eta_i(x)$ is the $i^{th}$-component of the unitary outer normal vector to $\partial\Omega$ in $x \in \partial\Omega$;

\item $u_0 \in \ W^{2 - \frac{2}{p}}_p(\Omega)$ with $p\neq 3$ and satisfying the compatibility condition
\\
$\displaystyle \sum\limits_{i=1}^n b_i \frac{\partial u_0}{\partial x_i} + b \ u_0 =0$ on $\partial \Omega$ when $p > 3$.

\end{enumerate}
Then, there exists a unique solution $u \in W^{2,1}_p(Q)$ of Problem~(\ref{P_Newmman});
moreover, there is a positive constant $C_p$  such that
the solution satisfies
\begin{equation}
\label{BasicParabolicEstimate}
\|u\|_{W^{2,1}_p(Q)} \leq C_{p} \left(\|f\|_{L^p(Q)} + \|u_0\|_{W^{2  - \frac{2}{p}}_p(\Omega)}\right).
\end{equation}

Such constant $C_{p}$ depends only on  $\Omega$, $T$, $p$, $r$, $s$, $\beta$, $\delta$
and on the norms $\|b_i\|_{C^2 (\bar{\Gamma})}$, $\|b\|_{C^2 (\bar{\Gamma})}$, $\|a_{ij}\|_{C(\bar{Q})}$,
$\|a_i\|_{L^r(Q)}$ and $\|a\|_{L^s(Q)}$.
Moreover, we may assume that the dependencies of $C_{p}$ on stated the norms are non decreasing.
\end{proposition}

\begin{remark}
The result set out in Proposition \ref{sol. Neumann} can be formulated for the parabolic problem with Dirichlet conditions (see Ladyzenskaya \cite[Theorem 9.1, pp.$341$]{Ladyzhenskaya}). In the problem with Dirichlet condition the compatibility condition in Proposition \ref{sol. Neumann}-($6$) can be replaced by $u_0=0$ on $\partial \Omega$ when $p > 3/2$. This way, all the results in this paper holds if we replaced the Neumann conditions by Dirichlet conditions.
\end{remark}

\section{An auxiliary problem}

In this section we will prove an auxiliary result  to be used in the proof of Theorem~\ref{Teorema1}.
To cope with difficulties with the signs of certain terms during the derivation of the estimates,
we firstly have to consider the following modified problem:
\begin{equation}
\label{P01}
\left\{
\begin{array}{lcl}
\displaystyle
\frac{\partial \hat{N}}{\partial t} = r_N - \mu_N \hat{N} - \beta_1 \hat{N} |\hat{A}| - \alpha_H\gamma_H \hat{N} |\hat{H}|, &\textup{in}& Q,
\vspace{0.2cm}
\\
\displaystyle
\frac{\partial \hat{A}}{\partial t} = \xi_A \Delta \hat{A} + r_A \hat{A}\left(1-\frac{\hat{A}}{k_A}\right)-(\mu_A+\epsilon_A)\hat{A} - \beta_3 \hat{N} \hat{A}  , &\textup{in}& Q,
\vspace{0.2cm}
\\
\displaystyle
\frac{\partial \hat{H}}{\partial t} = \xi_H \Delta \hat{H} + \nu \hat{A} - \tau_H \hat{H} - \gamma_H \hat{N} \hat{H} , &\textup{in}& Q,
\vspace{0.2cm}
\\
\displaystyle \frac{\partial \hat{A}}{\partial \eta} (\cdot) =\frac{\partial \hat{H}}{\partial \eta} (\cdot) =0, &\textup{on}& \Gamma,
\vspace{0.2cm}
\\
\displaystyle \hat{N}(\cdot,0) = N_0(\cdot), \hat{A}(\cdot,0) = A_{0}(\cdot), \hat{H}(\cdot,0) = H_0(\cdot), &\textup{in}& \Omega.
\end{array}
\right.
\end{equation}

Now we observe that, since the equation for $\hat{N}$ in this last problem is, for each $x \in \Omega$, an ordinary differential equation
which is linear in $\hat{N}$, we can find an explicit expression for it in terms of $|\hat{A}|$ and $|\hat{H}|$.
Using this observation, we introduce the operator $\Theta: L^{\infty}(Q) \times L^{\infty}(Q) \to L^{\infty}(Q)$, defined by
\begin{equation}
\label{P7}
\begin{array}{c}
\displaystyle	
\Theta(\phi, \varphi)(x,t) = \frac{N_0(x) + r_N \int_{0}^{t} e^{\mu_N s} e^{\alpha_H\gamma_H\int_{0}^{s}|\phi(x,\xi)| d\xi}e^{\beta_1\int_{0}^{s}|\varphi(x, \xi)| d\xi} ds}{e^{\mu_N t} e^{\alpha_H\gamma_H\int_{0}^{t}|\phi(x,\xi)| d\xi}e^{\beta_1\int_{0}^{t}|\varphi(x,\xi)| d\xi}},
\end{array}
\end{equation}

\noindent
where $0 \leq s \leq t \leq T$.

\begin{remark}
\label{obs1}
Thus,  $(\hat{N},\hat{A}, \hat{H})$ is solution of (\ref{P01}) if, and only if, $\hat{N} = \Theta (\hat{H}, \hat{A})$, $\hat{A}$ and $\hat{H}$ satisfies the following integro-differential system:
\begin{equation}
\label{P3}
\left\{
\begin{array}{lcl}
\displaystyle
\frac{\partial \hat{A}}{\partial t} = \xi_A \Delta \hat{A} + r_A \hat{A}\left(1-\frac{\hat{A}}{k_A}\right)-(\mu_A+\epsilon_A)\hat{A} - \beta_3 \Theta(\hat{A}, \hat{H}) \hat{A}  , &\textup{in}& Q,
\vspace{0.2cm}
\\
\displaystyle
\frac{\partial \hat{H}}{\partial t} = \xi_H \Delta \hat{H} + \nu \hat{A} - \tau_H \hat{H} - \gamma_H \Theta(\hat{A}, \hat{H}) \hat{H} , &\textup{in}& Q,
\vspace{0.2cm}
\\
\displaystyle \frac{\partial \hat{A}}{\partial \eta} (\cdot)= \frac{\partial \hat{H}}{\partial \eta} (\cdot) =0, &\textup{on}& \Gamma,
\vspace{0.2cm}
\\
\displaystyle \hat{A}(\cdot,0) = A_{0}(\cdot), \hat{H}(\cdot,0) = H_0(\cdot), &\textup{in}& \Omega.
\end{array}
\right.
\end{equation}
\end{remark}

\begin{remark}
\label{obs2}
Notice that, to guarantee that $(N, A, H)$, with $A= \hat{A}$, $H = \hat{H}$ and $N = \Theta(\hat{H}, \hat{A})$ is also a solution of system~(\ref{0riginalEquations}),
it is enough to prove that the solution $(\hat{A}, \hat{H})$ of Problem~(\ref{P3}) is nonnegative.
\end{remark}

For the Problem \ref{P3}, we have the following existence result:
\begin{proposition}\label{Prop1}
Assuming that the Technical Hypotheses~\ref{MainHypotheses} hold, there exists at least one nonnegative solution
$(\hat{A},\hat{H}) \in W^{2,1}_4(Q)\times W^{2,1}_4(Q)$ of Problem \eqref{P3}.
Moreover, $\hat{A} \leq k_A$ a.e. in $Q$ and
\begin{eqnarray*}
   ||\hat{A}||_{W^{2,1}_4(Q)} + ||\hat{H}||_{W^{2,1}_4(Q)} \leq C,
\end{eqnarray*}
where $C$ is a constant depending on $k_A$, $\nu$, $T$, $\Omega$, $||A_0||_{W^{\frac{3}{2}}_4(\Omega)}$ and $||H_0||_{W^{\frac{3}{2}}_4(\Omega)}$.
\end{proposition}

\begin{lemma}
\label{base}
Let $f:(0,T) \to \mathbb{R}$ differentiable such that $f(t) > 0$ and $f'(t) \ge 0$. If $g(t) = \frac{\int_{0}^{t} f(x) dx}{f(t)}$, then $g(t) \leq T$, for all $t \in (0, T)$.
\end{lemma}

\noindent
{\bf Proof:}
Since $f$ is continuous in $(0, T)$, it follows that
\begin{eqnarray*}
  g'(t)=\frac{f(t)^2 - f'(t) \int_{0}^{t}f(x)dx}{f(t)^2} = 1 - \frac{f'(t)}{f(t)} g(t).
\end{eqnarray*}

As $f(t) > 0$ we have $g(t) \ge 0$ and using the fact that $f'(t) \ge 0$ we obtain $\frac{f'(t)}{f(t)} g(t) \ge 0$. Therefore, $g'(t) \leq 1$, which suggests $g(t) \leq t$, for all $t \in (0, T)$. Thus, $g(t) \leq T$, as intended.
\hfill$\Box$

Since in the proof of existence of solutions of (\ref{P3})
the expression of $\Theta$ will play important roles,
we state in the following some of their properties:
\begin{lemma}
\label{PropertiesEtcFirst}
If $N_0 \in L^\infty(\Omega)$, then for any $\phi, \phi_1, \phi_2, \varphi, \varphi_1, \varphi_2  \in L^\infty(Q)$ and for almost everything $(x,t) \in Q$, there holds
\[
\begin{array}{ll}
\mbox{\bf (i)} & 0 \leq \Theta(\phi, \varphi)(x,t) \leq ||N_0||_{L^\infty(\Omega)} + r_N T;
\vspace{0.2cm}
\\
\mbox{\bf (ii)} & \|\Theta (\phi_1, \varphi_1) - \Theta (\phi_2, \varphi_2) \|_{L^\infty {(Q)}}  \leq C_1 \|(\phi_1, \varphi_1) - (\phi_2, \varphi_2)\|_{L^\infty(Q)}, \\
& where \ C_1 \ is \ a \ constant \ depending \ on \ r_N, \mu_N, \beta_1, \alpha_H, \gamma_H, T, ||\phi_1||_{L^\infty(Q)}, \\
& ||\phi_2||_{L^\infty(Q)}, ||\varphi_1||_{L^\infty(Q)}, ||\varphi_2||_{L^\infty(Q)} \ and \ ||N_0||_{L^\infty(\Omega)}.
\end{array}
\]
\end{lemma}

\noindent {\bf Proof (i):}
By expression (\ref{P7}) it is immediate that $\Theta(\phi, \varphi)(x,t) \ge 0$. To prove that $\Theta(\phi)(x,t) \leq ||N_0||_{L^\infty(\Omega)} + r_N T,$ we observe that
$$
\begin{array}{rcl}
\displaystyle
\Theta(\phi, \varphi)(x,t) = \frac{N_0(x) + r_N \int_{0}^{t} e^{\mu_N s} e^{\alpha_H\gamma_H\int_{0}^{s}|\phi(x,\xi)| d\xi}e^{\beta_1\int_{0}^{s}|\varphi(x, \xi)| d\xi} ds}{e^{\mu_N t} e^{\alpha_H\gamma_H\int_{0}^{t}|\phi(x,\xi)| d\xi}e^{\beta_1\int_{0}^{t}|\varphi(x,\xi)| d\xi}}   \\
\\
\displaystyle
\leq N_0(x) + r_N \frac{ \int_{0}^{t} e^{\mu_N s} e^{\alpha_H\gamma_H\int_{0}^{s}|\phi(x,\xi)| d\xi}e^{\beta_1\int_{0}^{s}|\varphi(x, \xi)| d\xi} ds}{e^{\mu_N t} e^{\alpha_H\gamma_H\int_{0}^{t}|\phi(x,\xi)| d\xi}e^{\beta_1\int_{0}^{t}|\varphi(x,\xi)| d\xi}} .
\end{array}
$$

Fixed $x \in \Omega$, we define
\begin{equation*}
\displaystyle
g(x,t) = \frac{\int_{0}^{t} e^{\mu_N s} e^{\alpha_H\gamma_H\int_{0}^{s}|\phi(x,\xi)| d\xi}e^{\beta_1\int_{0}^{s}|\varphi(x, \xi)| d\xi} ds}{e^{\mu_N t} e^{\alpha_H\gamma_H\int_{0}^{t}|\phi(x,\xi)| d\xi}e^{\beta_1\int_{0}^{t}|\varphi(x,\xi)| d\xi}},
\end{equation*}
and using the Lemma \ref{base} with $f(x,t) = e^{\mu_N t} e^{\alpha_H\gamma_H\int_{0}^{t}|\phi(x,\xi)| d\xi}e^{\beta_1\int_{0}^{t}|\varphi(x,\xi)| d\xi} $, it follow that
\begin{eqnarray*}
\Theta(\phi, \varphi)(x,t) \leq N_0(x) + r_N T \\
\leq ||N_0||_{L^\infty(\Omega)} + r_N T.
\end{eqnarray*}
{\bf Proof (ii):} We firstly need to observe that, due to the mean value inequality, given any $z_1, z_2 \in \R$,
there is $\theta = \theta(z_1, z_2)$ such that $e^{z_2} - e^{z_1} = e^{(1-\theta) z_1 + \theta z_2 } (z_2 - z_1)$;
in particular, for any $z_1, z_2 \leq 0$ we also have $(1-\theta) z_1 + \theta z_2  \leq 0$ and thus
\begin{equation}
\label{AlgebraicExponentialInequality}
|e^{z_2} - e^{z_1} | \leq |z_2 - z_1|, \quad \forall z_1, z_2 \leq 0 .
\end{equation}

Secondly, we note that
\begin{eqnarray*}
\displaystyle
\big|e^{\alpha_h\gamma_h \int_{0}^{t}|\phi_2(x,\xi)|d\xi}e^{\beta_1 \int_{0}^{t}|\varphi_2(x,\xi)|d\xi} - e^{\alpha_h\gamma_h \int_{0}^{t}|\phi_1(x,\xi)|d\xi}e^{\beta_1 \int_{0}^{t}|\varphi_1(x,\xi)|d\xi}\big| &\leq& \\
\displaystyle
e^{\alpha_h\gamma_h \int_{0}^{t}|\phi_2(x,\xi)|d\xi}\big|e^{\beta_1 \int_{0}^{t}|\varphi_2(x,\xi)|d\xi} - e^{\beta_1 \int_{0}^{t}|\varphi_1(x,\xi)|d\xi}\big| &+& \\
\displaystyle
e^{\beta_1 \int_{0}^{t}|\varphi_1(x,\xi)|d\xi} \big|e^{\alpha_h\gamma_h \int_{0}^{t}|\phi_2(x,\xi)|d\xi} - e^{\alpha_h\gamma_h \int_{0}^{t}|\phi_1(x,\xi)|d\xi}\big|.
\end{eqnarray*}

By the inequality (\ref{AlgebraicExponentialInequality}) and by $\phi_i, \varphi_i \in L^\infty(Q)$, $i = 1,2$, we obtain
\begin{equation}\label{i1}
\begin{array}{rcl}
\displaystyle
\big|e^{\alpha_h\gamma_h \int_{0}^{t}|\phi_2(x,\xi)|d\xi}e^{\beta_1 \int_{0}^{t}|\varphi_2(x,\xi)|d\xi} - e^{\alpha_h\gamma_h \int_{0}^{t}|\phi_1(x,\xi)|d\xi}e^{\beta_1 \int_{0}^{t}|\varphi_1(x,\xi)|d\xi}\big| &\leq& \\
\\
\displaystyle
e^{\alpha_h\gamma_h \int_{0}^{t}|\phi_2(x,\xi)|d\xi}\beta_1 \int_{0}^{t}|\varphi_2(x,\xi) - \varphi_1(x,\xi)|d\xi &+& \\
\\
\displaystyle
e^{\beta_1 \int_{0}^{t}|\varphi_1(x,\xi)|d\xi} \alpha_h\gamma_h \int_{0}^{t}|\phi_2(x,\xi) - \phi_1(x,\xi)|d\xi &\leq& \\
\\
\displaystyle
e^{\alpha_h\gamma_h \int_{0}^{t}|\phi_2(x,\xi)|d\xi}\beta_1 T ||\varphi_2 - \varphi_1||_{L^\infty(Q)} &+& \\
\\
\displaystyle
e^{\beta_1 \int_{0}^{t}|\varphi_1(x,\xi)|d\xi} \alpha_h\gamma_h T ||\phi_2 - \phi_1||_{L^\infty(Q)}.
\end{array}
\end{equation}

Thirdly, we observe that
\begin{eqnarray*}
\bigg|e^{\alpha_h\gamma_h \int_{0}^{t}|\phi_2(x,\xi)|d\xi}e^{\beta_1 \int_{0}^{t}|\varphi_2(x,\xi)|d\xi} \int_{0}^{t} e^{\mu_N s} e^{\alpha_H\gamma_H\int_{0}^{s}|\phi_1(x,\xi)| d\xi}e^{\beta_1\int_{0}^{s}|\varphi_1(x, \xi)| d\xi} ds &-& \\
\displaystyle
e^{\alpha_h\gamma_h \int_{0}^{t}|\phi_1(x,\xi)|d\xi}e^{\beta_1 \int_{0}^{t}|\varphi_1(x,\xi)|d\xi} \int_{0}^{t} e^{\mu_N s} e^{\alpha_H\gamma_H\int_{0}^{s}|\phi_2(x,\xi)| d\xi}e^{\beta_1\int_{0}^{s}|\varphi_2(x, \xi)| d\xi} ds \bigg| &\leq& \\
\displaystyle
e^{\alpha_h\gamma_h \int_{0}^{t}|\phi_2(x,\xi)|d\xi}e^{\beta_1 \int_{0}^{t}|\varphi_2(x,\xi)|d\xi} &\times& \\
\displaystyle
\int_{0}^{t} e^{\mu_N s}\big|e^{\alpha_H\gamma_H\int_{0}^{s}|\phi_1(x,\xi)| d\xi}e^{\beta_1\int_{0}^{s}|\varphi_1(x, \xi)| d\xi} - e^{\alpha_H\gamma_H\int_{0}^{s}|\phi_2(x,\xi)| d\xi}e^{\beta_1\int_{0}^{s}|\varphi_2(x, \xi)| d\xi}\big|ds &+& \\
\displaystyle
\int_{0}^{t} e^{\mu_N s} e^{\alpha_H\gamma_H\int_{0}^{s}|\phi_2(x,\xi)| d\xi}e^{\beta_1\int_{0}^{s}|\varphi_2(x, \xi)| d\xi} ds &\times& \\
\displaystyle
\big|e^{\alpha_h\gamma_h \int_{0}^{t}|\phi_2(x,\xi)|d\xi}e^{\beta_1 \int_{0}^{t}|\varphi_2(x,\xi)|d\xi} - e^{\alpha_h\gamma_h \int_{0}^{t}|\phi_1(x,\xi)|d\xi}e^{\beta_1 \int_{0}^{t}|\varphi_1(x,\xi)|d\xi}\big|.
\end{eqnarray*}

Study analogous to that done in (\ref{i1}), prove that
\begin{eqnarray*}
\bigg|e^{\alpha_h\gamma_h \int_{0}^{t}|\phi_2(x,\xi)|d\xi}e^{\beta_1 \int_{0}^{t}|\varphi_2(x,\xi)|d\xi} \int_{0}^{t} e^{\mu_N s} e^{\alpha_H\gamma_H\int_{0}^{s}|\phi_1(x,\xi)| d\xi}e^{\beta_1\int_{0}^{s}|\varphi_1(x, \xi)| d\xi} ds &-& \\
\displaystyle
e^{\alpha_h\gamma_h \int_{0}^{t}|\phi_1(x,\xi)|d\xi}e^{\beta_1 \int_{0}^{t}|\varphi_1(x,\xi)|d\xi} \int_{0}^{t} e^{\mu_N s} e^{\alpha_H\gamma_H\int_{0}^{s}|\phi_2(x,\xi)| d\xi}e^{\beta_1\int_{0}^{s}|\varphi_2(x, \xi)| d\xi} ds \bigg| &\leq& \\
\displaystyle
e^{\alpha_h\gamma_h \int_{0}^{t}|\phi_2(x,\xi)|d\xi}e^{\beta_1 \int_{0}^{t}|\varphi_2(x,\xi)|d\xi} &\times& \\
\displaystyle
e^{\mu_N T} T \bigg(e^{\alpha_h\gamma_h T ||\phi_2||_{L^\infty(Q)}} \beta_1 T ||\varphi_2 - \varphi_1||_{L^\infty(Q)} +
e^{\beta_1 T ||\varphi_1||_{L^\infty(Q)}} \alpha_h\gamma_h T ||\phi_2 - \phi_1||_{L^\infty(Q)} \bigg) &+& \\
\displaystyle
T e^{\mu_N T} e^{\alpha_H\gamma_H T ||\phi_2||_{L^\infty(Q)}}e^{\beta_1 T ||\varphi_2||_{L^\infty(Q)}} &\times& \\
\displaystyle
\bigg(e^{\alpha_h\gamma_h \int_{0}^{t}|\phi_2(x,\xi)|d\xi}\beta_1 T ||\varphi_2 - \varphi_1||_{L^\infty(Q)} +
e^{\beta_1 \int_{0}^{t}|\varphi_1(x,\xi)|d\xi} \alpha_h\gamma_h T ||\phi_2 - \phi_1||_{L^\infty(Q)} \bigg).
\end{eqnarray*}

Finally, the expression in (\ref{P7}) suggests
\begin{equation}
\label{ineq1}
\begin{array}{rcl}
\displaystyle
|\Theta(\phi_1, \varphi_1)(x,t) - \Theta(\phi_2, \varphi_2)(x,t)| & \leq & \\
\\
\displaystyle
\frac{1}{e^{\mu_N t}e^{\alpha_h\gamma_h \int_{0}^{t}|\phi_1(x,\xi)|d\xi}e^{\beta_1 \int_{0}^{t}|\varphi_1(x,\xi)|d\xi}e^{\alpha_h\gamma_h \int_{0}^{t}|\phi_2(x,\xi)|d\xi}e^{\beta_1 \int_{0}^{t}|\varphi_2(x,\xi)|d\xi}} &\times& \\
\\
\displaystyle
\bigg(N_0(x)e^{\mu_n t}\bigg|e^{\alpha_h\gamma_h \int_{0}^{t}|\phi_2(x,\xi)|d\xi}e^{\beta_1 \int_{0}^{t}|\varphi_2(x,\xi)|d\xi} - e^{\alpha_h\gamma_h \int_{0}^{t}|\phi_1(x,\xi)|d\xi}e^{\beta_1 \int_{0}^{t}|\varphi_1(x,\xi)|d\xi}\bigg| &+& \\
\\
\displaystyle
r_N e^{\mu_N t} \bigg|e^{\alpha_h\gamma_h \int_{0}^{t}|\phi_2(x,\xi)|d\xi}e^{\beta_1 \int_{0}^{t}|\varphi_2(x,\xi)|d\xi} \int_{0}^{t} e^{\mu_N s} e^{\alpha_H\gamma_H\int_{0}^{s}|\phi_1(x,\xi)| d\xi}e^{\beta_1\int_{0}^{s}|\varphi_1(x, \xi)| d\xi} ds &-&
\\
\displaystyle
e^{\alpha_h\gamma_h \int_{0}^{t}|\phi_1(x,\xi)|d\xi}e^{\beta_1 \int_{0}^{t}|\varphi_1(x,\xi)|d\xi} \int_{0}^{t} e^{\mu_N s} e^{\alpha_H\gamma_H\int_{0}^{s}|\phi_2(x,\xi)| d\xi}e^{\beta_1\int_{0}^{s}|\varphi_2(x, \xi)| d\xi} ds \bigg|\bigg)
\end{array}
\end{equation}
and using the estimates obtained in (\ref{i1}) and (\ref{ineq1}) and making the possible simplifications, we obtain
$$
\begin{array}{rcl}
\displaystyle
|\Theta(\phi_1, \varphi_1)(x,t) - \Theta(\phi_2, \varphi_2)(x,t)| & \leq & \\
\\
\displaystyle
||N_0||_{L^\infty(\Omega)} \big( \beta_1 T ||\varphi_2 - \varphi_1||_{L^\infty(Q)} + \alpha_h\gamma_h T ||\phi_2 - \phi_1||_{L^\infty(Q)} \big) &+& \\
\\
\displaystyle
r_N \bigg(e^{\mu_N T} T^2 e^{\alpha_h\gamma_h ||\phi_2||_{L^\infty(Q)}}\beta_1 ||\varphi_2 - \varphi_1||_{L^\infty(Q)} &+& \\
\\
\displaystyle
e^{\mu_N T} T^2 e^{\beta_1 T ||\varphi_1||_{L^\infty(Q)}} \alpha_h\gamma_h ||\phi_2 - \phi_1||_{L^\infty(Q)} &+& \\
\\
\displaystyle
\beta_1 T^2 e^{\mu_N T} e^{\alpha_H\gamma_H T ||\phi_2||_{L^\infty(Q)}} e^{\beta_1 T ||\varphi_2||_{L^\infty(Q)}} ||\varphi_2 - \varphi_1||_{L^\infty(Q)} &+& \\
\\
\displaystyle
\alpha_h\gamma_h  T^2 e^{\mu_N T} e^{\alpha_H\gamma_H T ||\phi_2||_{L^\infty(Q)}} e^{\beta_1 T ||\varphi_2||_{L^\infty(Q)}} ||\phi_2 - \phi_1||_{L^\infty(Q)}\bigg),
\end{array}
$$
for almost everything $(x,t) \in Q$, i.e.,
\begin{eqnarray*}
\displaystyle
||\Theta(\phi_1, \varphi_1) - \Theta(\phi_2, \varphi_2)||_{L^\infty(Q)} \leq C_1 \big(||\phi_2 - \phi_1||_{L^\infty(Q)} + ||\varphi_2 - \varphi_1||_{L^\infty(Q)} \big).
\end{eqnarray*}
\hfill$\Box$

\subsection{Proof of Proposition \ref{Prop1}}

To not overburden the notation, in this subsection we denote $(A,H)$ as a generic solution of the equations that follows.

To get a solution of problem \eqref{P3}, we will apply the Leray-Schauder fixed point theorem to the mapping
$\Psi$ defined as follows:
\begin{equation}
\label{oper}
\begin{array}{rccl}
\Psi: & [0,1]\times L^\infty(Q) \times L^\infty(Q)  & \rightarrow & L^\infty(Q) \times L^\infty(Q)
\\
&(l, \phi, \varphi) & \mapsto &	(A,H) ,
\end{array}
\end{equation}

\noindent
where $(A,H)$ is the unique solution of
\begin{equation}
\label{P66}
\left\{
\begin{array}{lcl}
\displaystyle
\frac{\partial A}{\partial t} = \xi_A \Delta A + r_A A\left(1-\frac{A}{k_A}\right)-(\mu_A+\epsilon_A)A - l \beta_3 \Theta(\phi, \varphi) A  , &\textup{in}& Q,
\vspace{0.2cm}
\\
\displaystyle
\frac{\partial H}{\partial t} = \xi_H \Delta H + \nu A - \tau_H H - l \gamma_H \Theta(\phi, \varphi) H , &\textup{in}& Q,
\vspace{0.2cm}
\\
\displaystyle \frac{\partial A}{\partial \eta} (\cdot) =\frac{\partial H}{\partial \eta} (\cdot) =0, &\textup{on}& \Gamma,
\vspace{0.2cm}
\\
\displaystyle A(\cdot,0) = A_{0}(\cdot), H(\cdot,0) = H_0(\cdot), &\textup{in}& \Omega,
\end{array}
\right.
\end{equation}
with $\Theta(\phi, \varphi)$ given by \eqref{P7}.

We affirm that the operator defined in (\ref{oper}) is well defined if the problems
\begin{equation}
\label{separado1}
\left\{
\begin{array}{lcl}
\displaystyle
\frac{\partial A}{\partial t} = \xi_A \Delta A + r_A A\left(1-\frac{A}{k_A}\right)-(\mu_A+\epsilon_A)A - l\beta_3 \Theta(\phi, \varphi) A  , &\textup{in}& Q,
\vspace{0.2cm}
\\
\displaystyle \frac{\partial A}{\partial \eta} (\cdot) =0, &\textup{on}& \Gamma,
\vspace{0.2cm}
\\
\displaystyle A(\cdot,0) = A_{0}(\cdot), &\textup{in}& \Omega
\end{array}
\right.
\end{equation}
and
\begin{equation}
\label{separado2}
\left\{
\begin{array}{lcl}
\displaystyle
\frac{\partial H}{\partial t} = \xi_H \Delta H + \nu A - \tau_H H - l \gamma_H \Theta(\phi, \varphi) H , &\textup{in}& Q,
\vspace{0.2cm}
\\
\displaystyle \frac{\partial H}{\partial \eta} (\cdot) =0, &\textup{on}& \Gamma,
\vspace{0.2cm}
\\
\displaystyle H(\cdot,0) = H_0(\cdot), &\textup{in}& \Omega,
\end{array}
\right.
\end{equation}
have a unique solution.

\subsection{Existence and uniqueness of solution to problems (\ref{separado1}) and (\ref{separado2})}
\label{sec32}

\noindent\textbf{Step 1:} To determine a solution for the problem (\ref{separado1}), start by studying the following modified problem:
\begin{equation}
\label{separado31}
\left\{
\begin{array}{lcl}
\displaystyle
\frac{\partial A}{\partial t} = \xi_A \Delta A + r_A A \bigg(1 - \frac{A^+}{k_A}\bigg) - (\mu_A + \epsilon_A)A - l \beta_3 \Theta(\phi, \varphi) A, &\textup{in}& Q,
\vspace{0.2cm}
\\
\displaystyle \frac{\partial A}{\partial \eta} (\cdot) =0, &\textup{on}& \Gamma,
\vspace{0.2cm}
\\
\displaystyle A(\cdot,0) = A_{0}(\cdot), &\textup{in}& \Omega.
\end{array}
\right.
\end{equation}

In order to prove that problem (\ref{separado31}) admits solution, we define the operator
\begin{equation}
\label{oper1}
\begin{array}{rcl}
G: [0,1] \times L^\infty(Q) &\rightarrow& L^\infty(Q) \\
\\
(\lambda, g) &\mapsto&  A,
\end{array}
\end{equation}
where $A$ is the unique solution to the problem
\begin{equation}
\label{separado33}
\left\{
\begin{array}{lcl}
\displaystyle
\frac{\partial A}{\partial t} = \xi_A \Delta A + \lambda r_A A \bigg(1 - \frac{g^{+}}{k_A}\bigg) - (\mu_A+\epsilon_A)A - l \beta_3 \Theta(\phi, \varphi) A, &\textup{in}& Q,
\vspace{0.2cm}
\\
\displaystyle \frac{\partial A}{\partial \eta} (\cdot) =0, &\textup{on}& \Gamma,
\vspace{0.2cm}
\\
\displaystyle A(\cdot,0) = A_{0}(\cdot), &\textup{in}& \Omega.
\end{array}
\right.
\end{equation}

\vspace{0.1cm}
To ensure that we can apply the Leray-Schauder fixed point theorem, we present next a sequence of lemmas:

\begin{lemma}
\label{lemaAi}
Suppose $N_0 \in L^\infty(\Omega)$ and $A_0 \in W^{\frac{3}{2}}_4(\Omega)$. Then the mapping $G :[0,1] \times L^\infty(Q) \rightarrow L^\infty(Q)$ is well defined.
\end{lemma}

\noindent{\bf Proof:}
Note that the coefficients of the problem \eqref{separado33} satisfy the hypotheses of the Proposition \ref{sol. Neumann}. For example, it is immediate that $\lambda r_A - \lambda \frac{r_A}{k_A}g^+ - (\mu_A+\epsilon_A) - l \beta_3 \Theta(\phi, \varphi) \in L^4(Q)$, because $g \in L^\infty(Q)$ and by Lemma \ref{PropertiesEtcFirst}, $\Theta(\phi, \varphi) \in L^\infty(Q)$. Thus, we conclude that there is a unique solution $A \in W^{2,1}_4(Q)$ of problem \eqref{separado33}. Moreover, $A$ satisfies the following estimate:
\begin{equation}\label{AA}
  ||A||_{W^{2,1}_4(Q)} \leq C_p ||A_0||_{W^{\frac{3}{2}}_4(\Omega)}.
\end{equation}

Finally, from Lemma \ref{icontLp01}, we have $W^{2,1}_4(Q) \hookrightarrow L^\infty(Q)$, and we conclude that the operator $G$ in well defined.
\hfill$\Box$
\\

\begin{lemma}
\label{nãonegativa1}
Suppose $A$ a solution of (\ref{separado33}) and $A_0 \ge 0$ a.e. in $\Omega$, then $A \ge 0$ a.e. in $Q$.
\end{lemma}

\noindent {\bf Proof:}
Multiplying the first equation in $(\ref{separado33})$ by $A^-$ and integrating into $\Omega$, we get
\begin{eqnarray*}
-\frac{1}{2} \frac{d}{dt} \int_{\Omega} (A^{-})^2 dx  &=&  \xi_A \int_{\Omega} |\nabla A^-|^2 dx - \lambda r_A \int_{\Omega} (A^-)^2 dx + (\mu_A+\epsilon_A) \int_{\Omega}( A^-)^2 dx \\
&+&  \lambda  \frac{r_A}{k_A} \int_{\Omega} g_+ (A^-)^2 dx +
l \beta_3 \int_{\Omega}  \Theta(\phi, \varphi) (A^-)^2 dx.
\end{eqnarray*}

Thus,
\begin{eqnarray*}
\frac{d}{dt} \int_{\Omega} (A^{-})^2 dx \leq 2 r_A \int_{\Omega} (A^-)^2 dx
\end{eqnarray*}
and using the Gronwall's inequality and the fact that $A_0 \ge 0$ a.e. in $\Omega$, we obtain
\begin{eqnarray*}
\int_{\Omega} (A^{-})^2 dx \leq e^{2 r_A T} \int_{\Omega}
({A_0}^-)^2 dx = 0,
\end{eqnarray*}
that is, $||A^-(\cdot, t)||_{L^2(\Omega)} = 0$ for all $t \in
(0,T)$, where we conclude that $A^{-} = 0$ a.e. in $Q$ and therefore $A
\ge 0$ a.e. in $Q$.
\hfill$\Box$
\\

\begin{lemma}
\label{con1}
For each fixed $\lambda\in[0,1]$, the mapping $G(\lambda, \cdot): L^\infty(Q) \rightarrow L^\infty(Q)$ is compact, i.e.,
it is continuous and maps bounded sets into relatively compacts sets.
\end{lemma}

\noindent{\bf Proof:}
The functions $G(\lambda, g_1)= A_1$ and $G(\lambda, g_2) = A_2$ satisfy the problem
\begin{equation*}
\label{Aij}
\left\{
\begin{array}{lcl}
\displaystyle \frac{\partial A_i}{\partial t} = \xi_A \Delta
A_i + \lambda r_A A_i\bigg(1 - \frac{g^{+}_i}{k_A} \bigg) - (\mu_A + \epsilon_A) A_i - l\beta_3\Theta(\phi, \varphi) A_i , &\textup{in}& Q,
\vspace{0.2cm}
\\
\displaystyle \frac{\partial A_i}{\partial \eta} (\cdot) =0, &\textup{on}& \Gamma,
\vspace{0.2cm}
\\
\displaystyle A_i(\cdot,0) = A_0(\cdot), &\textup{in}& \Omega,
\end{array}
\right.
\end{equation*}
with $i=1,2$; letting $\tilde{A} := A_1 - A_2$, we have
\begin{equation}
\label{Ai1}
\left\{
\begin{array}{lcl}
\displaystyle \frac{\partial \tilde{A}}{\partial t} = \xi_A \Delta
\tilde{A} + \bigg(\lambda r_A - (\mu_A + \epsilon_A) - l \beta_3
\Theta(\phi, \varphi) - \lambda \frac{r_A}{k_A}{g_2}^+ \bigg)
\tilde{A} - 
\vspace{0.2cm}
\\
\displaystyle \lambda \frac{r_A}{k_A} A_1 ({g_1}^+ - {g_2}^+),
&\textup{in}& Q, 
\vspace{0.2cm}
\\
\displaystyle \frac{\partial \tilde{A}}{\partial \eta} (\cdot) =0, &\textup{on}& \Gamma,
\vspace{0.2cm}
\\
\displaystyle \tilde{A}(\cdot,0) = \tilde{A}_0(\cdot) = 0, &\textup{in}& \Omega.
\end{array}
\right.
\end{equation}

Since $g_i, A_i \in L^\infty(Q)$, $i=1,2$, from Proposition \ref{sol. Neumann}, we get
\begin{eqnarray*}
  ||\tilde{A}||_{W_{4}^{2, 1}(Q)} \leq C_p \lambda \frac{r_A}{k_A} || A_1 ({g_1}^+ - {g_2}^+) ||_{L^4(Q)} \\
  \leq C_p \frac{r_A}{k_A} T^{\frac{1}{4}} |\Omega|^{\frac{1}{4}} ||A_1||_{L^\infty(Q)} ||g_1 - g_2||_{L^\infty(Q)}.
\end{eqnarray*}

Then, by Lemma \ref{icontLp01}, we finally have
\begin{eqnarray*}
  ||G(\lambda, g_1) - G(\lambda, g_2)||_{L^\infty(Q)} \leq C ||g_1 - g_2||_{L^\infty(Q)},
\end{eqnarray*}
where $C$ depends on $C_p$, $r_A$, $k_A$, $T$, $\Omega$, $||A_1||_{L^\infty(Q)}$ and the immersion constant.

To show that $ G (\lambda, \cdot) $ is compact, we use the fact that immersion
$W^{2,1}_4(Q) \hookrightarrow L^\infty(Q)$ is compact and that
$G(\lambda, \cdot)$ is the composition between the inclusion operator and the
solution operator, i.e., $G(\lambda, \cdot): L^\infty(Q)
\rightarrow W^{2,1}_4(Q) \rightarrow L^\infty(Q)$.
\hfill$\Box$
\\

\begin{lemma}
\label{unifcon1}
Given a bounded subset $B\subset\L^\infty(Q)$, for each $g \in B$, the mapping $G(\cdot, g): [0, 1] \rightarrow L^\infty(Q)$ is continuous, uniformly with respect to $B$.
\end{lemma}

\noindent{\bf Proof:}
Since $B \in L^\infty(Q)$ is bounded, there is $r_B \ge 0$ such that, for any $g \in B$, we have $||g||_{L^\infty(Q)} \leq r_B$. Now, let us fix $g \in L^\infty(Q)$ and consider $\lambda_1, \lambda_2 \in [0,1]$ and denote $G(\lambda_1, g) = A_1$, $G(\lambda_2, g) = A_2$ and $\tilde{A} = A_1 - A_2$. Then, $\tilde{A}$ satisfies
\begin{equation}
\label{B}
\left\{
\begin{array}{lcl}
\displaystyle \frac{\partial \tilde{A}}{\partial t} = \xi_A \Delta
\tilde{A} - \bigg(r_A \lambda_2 - (\mu_A + \epsilon_A) - l \beta_3
\Theta(\phi, \varphi) - \frac{r_A}{k_A} \lambda_2 g^+ \bigg)
\tilde{A} - \vspace{0.2cm}
\\
\displaystyle \bigg(r_A A_1 + \frac{r_A}{k_A} g^+ A_1\bigg) (\lambda_1 - \lambda_2),
&\textup{in}& Q, \vspace{0.2cm}
\\
\displaystyle \frac{\partial \tilde{A}}{\partial \eta} = 0, &\textup{on}& \Gamma,
\vspace{0.2cm}
\\
\displaystyle \tilde{A}(\cdot, 0) = \tilde{A}_0(\cdot) = 0, &\textup{in}& \Omega.
\end{array}
\right.
\end{equation}

Since $g, A_i \in L^\infty(Q)$, $i=1,2$, from Proposition \ref{sol. Neumann}, we get
\begin{eqnarray*}
  ||\tilde{A}||_{W_{4}^{2,1}(Q)} \leq C_p |\lambda_1 - \lambda_2| \bigg|\bigg|\bigg(r_A A_1 + \frac{r_A}{k_A} g^+ A_1\bigg)\bigg|\bigg|_{L^4(Q)} \\
  \leq C_p \bigg(r_A + \frac{r_A}{k_A} r_B \bigg) T^{\frac{1}{4}} |\Omega|^{\frac{1}{4}} ||A_1||_{L^\infty(Q)} |\lambda_1 - \lambda_2|.
\end{eqnarray*}

Then, by Lemma \ref{icontLp01}, we finally have
\begin{eqnarray*}
  ||G(\lambda_1,g) - G(\lambda_2,g)||_{L^\infty(Q)} \leq C |\lambda_1 - \lambda_2|,
\end{eqnarray*}
where $C$ depends on $C_p$, $r_A$, $k_A$, $r_B$, $T$, $\Omega$, $||A_1||_{L^\infty(Q)}$ and the immersion constant.
\hfill$\Box$
\\

\begin{lemma}
\label{estimativa1}
Suppose $A_0 \leq k_A$ a.e. in $\Omega$, then there exists a number $\rho_1 > 0$ such that, for any $\lambda \in [0,1]$ and any possible
fixed point $A \in L^\infty(Q)$ of $G(\lambda, \cdot)$, there holds $\|A\|_{L^\infty(Q)} < \rho_1$.
\end{lemma}

\noindent {\bf Proof:} Let $A \in L^\infty(Q)$ be a possible fixed point of $G(l,\cdot)$ associated to a given $\lambda \in [0,1]$; then $A$ satisfies
\begin{equation}\label{rhox}
\left\{
\begin{array}{lcl}
\displaystyle \frac{\partial A}{\partial t} = \xi_A \Delta A + \lambda r_A A
\bigg(1 - \frac{A}{k_A}\bigg) - (\mu_A + \epsilon_A)A - l
\beta_3 \Theta(\phi, \varphi) A, &\textup{in}& Q, \vspace{0.2cm}
\\
\displaystyle \frac{\partial A}{\partial \eta} (\cdot) =0, &\textup{on}& \Gamma,
\vspace{0.2cm}
\\
\displaystyle A(\cdot,0) = A_{0}(\cdot), &\textup{in}& \Omega,
\end{array}
\right.
\end{equation}
since by the Lemma \ref{nãonegativa1}, $A = A^+$.

We observe that the first equation in (\ref{rhox}) can be rewritten as
\begin{eqnarray*}
\frac{\partial}{\partial t} (A - k_A) &=& \xi_A \Delta (A - k_A) - \lambda \frac{r_A}{k_A} A (A - k_A) \\
&-& (\mu_A+\epsilon_A)A - \beta_3 \Theta(\phi, \varphi) A.
\end{eqnarray*}

Now, multiplying by $(A - k_A)^+$ and integrating in $\Omega$, we obtain
\begin{eqnarray*}
\frac{1}{2} \frac{d}{dt} \int_{\Omega} \big((A - k_A)^+\big)^2 dx &=& - \xi_A \int_{\Omega} |\nabla(A - k_A)^{+}|^2 dx - \lambda \frac{r_A}{k_A} \int_{\Omega} A \big((A - k_A)^+\big)^2 dx \\
&-& (\mu_A+\epsilon_A) \int_{\Omega} A (A - k_A)^+ dx \\
&-& \beta_3 \int_{\Omega} \Theta(\phi, \varphi) A (A - k_A)^+ dx \\
&\leq& 0.
\end{eqnarray*}

Thus, using the Gronwall's inequality and the fact that $A_0 \leq k_A$ a.e. in $\Omega$, it follow that
\begin{eqnarray*}
\int_{\Omega} \big((A - k_A)^+\big)^2 dx
&\leq & \int_{\Omega} \big((A_0 - k_A)^+\big)^2 dx = 0,
\end{eqnarray*}
that is, $||(A(\cdot, t) - k_A)^+||_{L^2(\Omega)} = 0$ for all $t \in (0,T)$ and therefore $(A - k_A)^+ = 0$ a.e. in $Q$, and we conclude that $A \leq k_A$ a.e. in $Q$.

Thus, $||A||_{L^\infty(Q)} \leq k_A$, which suggests $\rho_1 = k_A + 1$.
\hfill$\Box$
\\

\begin{lemma}
\label{fix1}
The mapping $G(0,\cdot): L^\infty(Q) \rightarrow L^\infty(Q)$ has a unique fixed point.
\end{lemma}

\noindent {\bf Proof:}
Indeed, letting $\lambda =0$ in\eqref{separado33}, $A$ is a fixed point of $G(0, \cdot)$ if, and only if, $A$ is the unique solution to the problem
\begin{equation*}
\label{separado3}
\left\{
\begin{array}{lcl}
\displaystyle
\frac{\partial A}{\partial t} = \xi_A \Delta A - \big( (\mu_A+\epsilon_A) + l\beta_3\Theta(\phi, \varphi)\big)A, &\textup{in}& Q,
\vspace{0.2cm}
\\
\displaystyle \frac{\partial A}{\partial \eta} (\cdot) =0, &\textup{on}& \Gamma,
\vspace{0.2cm}
\\
\displaystyle A(\cdot,0) = A_{0}(\cdot), &\textup{in}& \Omega.
\end{array}
\right.
\end{equation*}
By Proposition {\ref{sol. Neumann}}, we have the existence of a unique solution $A \in W^{2,1}_{4}(Q)\hookrightarrow L^{\infty}(Q)$
of this problem; therefore $G(0, \cdot )$ has a unique fixed point in $L^\infty(Q)$.
\hfill$\Box$
\\

\begin{proposition}
\label{existenciaA}
There is a nonnegative solution $A \in W^{2,1}_4(Q)$ of the problem (\ref{separado1}).
\end{proposition}

\noindent {\bf Proof:}
From Lemmas \ref{lemaAi}, \ref{con1}, \ref{unifcon1}, \ref{estimativa1} and \ref{fix1}, we conclude that the mapping $G: [0,1] \times L^\infty(Q) \rightarrow L^\infty(Q)$ satisfies the hypotheses of the Leray-Schauder's fixed point theorem (see Friedman \cite[pp.~189, Theorem 3]{Friedman}). Thus, there exists $A \in L^\infty(Q)$ such that $G(1, A) = A$. Moreover, by Lemmas \ref{lemaAi} and \ref{nãonegativa1}, $A \in W^{2,1}_4(Q)$ is nonnegative and $A$ is the required solution of (\ref{separado1}).
\hfill$\Box$
\\

\begin{lemma}
\label{nãonegativa3}
Suppose $A$ a solution of (\ref{separado1}) and $A_0 \leq k_A$ a.e. in $\Omega$, then $A \leq k_A$ a.e. in $Q$.
\end{lemma}
\noindent {\bf Proof:}
Analogous to Lemma \ref{estimativa1}.
\hfill$\Box$
\\

\begin{proposition}
\label{unicidadeA}
The solution $A$ of the problem (\ref{separado1}) is unique.
\end{proposition}

\noindent {\bf Proof:}
Let $A_1$ and $A_2$ be solutions to the problem (\ref{separado1}); if $\tilde{A} = A_1 - A_2$, then $\tilde{A}$ satisfies the o problem
\begin{equation}
\label{unicidade11}
\left\{
\begin{array}{lcl}
\displaystyle
\frac{\partial \tilde{A}}{\partial t} = \xi_A \Delta \tilde{A} + r_A \tilde{A} - \frac{r_A}{k_A}(A_1 + A_2) \tilde{A} - (\mu_A+\epsilon_A) \tilde{A} - l\beta_3 \Theta(\phi, \varphi) \tilde{A}  , &\textup{in}& Q,
\vspace{0.2cm}
\\
\displaystyle \frac{\partial \tilde{A}}{\partial \eta} (\cdot) =0, &\textup{on}& \Gamma,
\vspace{0.2cm}
\\
\displaystyle \tilde{A}(\cdot,0) = \tilde{A}_0(\cdot) = 0, &\textup{in}& \Omega.
\end{array}
\right.
\end{equation}
Multiplying the first equation of (\ref{unicidade11}) by $\tilde{A}$ and integrating in $\Omega$, we obtain
\begin{eqnarray*}
  \frac{1}{2}\frac{d}{dt} \int_{\Omega} \tilde{A}^2 dx = - \xi_A \int_{\Omega} |\nabla \tilde{A}|^2 dx + r_A \int_{\Omega} \tilde{A}^2 dx - \frac{r_A}{k_A} \int_{\Omega} (A_1 + A_2) \tilde{A}^2 dx \\
  - (\mu_A+\epsilon_A) \int_{\Omega} \tilde{A}^2 dx - l\beta_3 \int_{\Omega} \Theta(\phi, \varphi) \tilde{A}^2 dx.
\end{eqnarray*}

Therefore,
\begin{eqnarray*}
  \frac{d}{dt} \int_{\Omega} \tilde{A}^2 dx &\leq& 2 r_A \int_{\Omega} \tilde{A}^2 dx,
\end{eqnarray*}
and by the Gronwall's inequality, we obtain
\begin{eqnarray*}
  \int_{\Omega} \tilde{A}^2 dx &\leq& e^{2r_A T} \int_{\Omega} (\tilde{A}_0)^2 dx = 0,
\end{eqnarray*}
that is, $||\tilde{A}(\cdot,t)||_{L^2(\Omega)} = 0$ for all $t \in (0, T)$, which suggests $\tilde{A} = 0$ a.e. in $Q$ and therefore $A_1 = A_2$ a.e. in $Q$.
\hfill$\Box$
\\

\noindent\textbf{Step 2:} The existence and uniqueness of solution for the problem (\ref{separado2}) is given by combining Propositions \ref{existenciaA} and \ref{unicidadeA} and the following result:

\begin{proposition}
\label{existenciaH}
Suppose $N_0 \in L^{\infty}(\Omega)$ and $H_0 \in W_{4}^{\frac{3}{2}}(\Omega)$. Then the problem (\ref{separado2}) has a unique solution.
\end{proposition}

\noindent {\bf Proof:}
It's immediate that the coefficients of the problem \eqref{separado2} satisfy the hypotheses of Proposition \ref{sol. Neumann}, and we conclude that there is a unique solution $H \in W^{2,1}_4(Q)$ of the problem $(\ref{separado2})$. Moreover, $H$ satisfies the following estimate:
\begin{equation}\label{HH}
  ||H||_{W^{2,1}_4(Q)} \leq C_p \bigg( \nu k_A T^{\frac{1}{4}}|\Omega|^{\frac{1}{4}} + ||H_0||_{W^{\frac{3}{2}}_4(\Omega)}\bigg)
\end{equation}
\hfill$\Box$
\\

\begin{lemma}
\label{nãonegativa2} Suppose $H$ a solution do problem (\ref{separado2}) and $H_0 \ge 0$ a.e. in $\Omega$, then $H \ge 0$ a.e. in $Q$.
\end{lemma}

\noindent {\bf Proof:}
Multiplying the first equation in (\ref{separado2}) by $H^-$ and integrating in $\Omega$, we get
\begin{eqnarray*}
-\frac{1}{2} \frac{d}{dt} \int_{\Omega} (H^-)^2 dx &=& \xi_H \int_{\Omega} |\nabla H^-|^2 dx + \nu \int_{\Omega} A H^- dx \\
&+& \tau_H \int_{\Omega} (H^-)^2 + l \gamma_H \int_{\Omega} \Theta(\phi, \varphi) (H^-)^2 dx.
\end{eqnarray*}

Using the fact that $A \ge 0$, we obtain
\begin{eqnarray*}
\frac{d}{dt} \int_{\Omega} (H_-)^2 dx \leq 0.
\end{eqnarray*}

Therefore, using the Gronwall's inequality and the fact that $H_0 \ge 0$ a.e. on $\Omega$, it follow that $||H^-(\cdot, t)||_{L^2(\Omega)} = 0$ for all $t \in (0,T)$, that is, $H^- = 0$ a.e. in $Q$, which suggests $H \ge 0$ a.e. in $Q$.
\hfill$\Box$
\\

\subsection{Continuation of Proof of Proposition \ref{Prop1}}
To guarantee that we can apply the Leray-Schauder fixed point theorem, we present next a sequence of lemmas.

\begin{lemma}
\label{bemdefinido}
The mapping $\Psi :[0,1]\times L^\infty(Q) \times L^\infty(Q)  \rightarrow L^\infty(Q) \times L^\infty(Q)$ is well defined.
\end{lemma}

\noindent{\bf Proof:}
Just combine the Propositions \ref{existenciaA}, \ref{unicidadeA}  and \ref{existenciaH}.
\hfill$\Box$
\\

\begin{lemma}
\label{con}
For each fixed $l \in[0,1]$, the mapping $\Psi(l, \cdot, \cdot): L^{\infty}(Q) \times L^{\infty}(Q) \rightarrow L^{\infty}(Q) \times L^{\infty}(Q)$ is compact, i.e., it is continuous and maps bounded sets into relatively compacts sets.
\end{lemma}

\noindent {\bf Proof:}
Consider $(\phi_1, \varphi_1), (\phi_2, \varphi_2) \in L^{\infty}(Q) \times L^{\infty}(Q)$ such that $\Psi(l, \phi_1, \varphi_1) = (A_1, H_1), \Psi(l, \phi_2, \varphi_2) = (A_2, H_2)$; if $\tilde{A} = A_1 - A_2$ and $\tilde{H} = H_1 - H_2$, we have that $\tilde{A}$ and $\tilde{H}$ satisfy the following problems, respectively:
\begin{equation}
\label{separado1i}
\left\{
\begin{array}{lcl}
\displaystyle
\frac{\partial \tilde{A}}{\partial t} = \xi_A \Delta \tilde{A} + \bigg(r_A - (\mu_A+\epsilon_A) - l \frac{r_A}{k_A}(A_1 + A_2) - l\beta_3\Theta(\phi_1, \varphi_1)\bigg) \tilde{A} -  
\vspace{0.2cm}
\\
\displaystyle
l\beta_3\big[\Theta(\phi_1, \varphi_1) - \Theta(\phi_1, \varphi_1)\big] A_2  , &\textup{in}& Q,
\vspace{0.2cm}
\\
\displaystyle \frac{\partial \tilde{A}}{\partial \eta} (\cdot) =0, &\textup{on}& \Gamma,
\vspace{0.2cm}
\\
\displaystyle \tilde{A}(\cdot,0) = \tilde{A}_{0}(\cdot) = 0, &\textup{in}& \Omega,
\end{array}
\right.
\end{equation}

\begin{equation}
\label{separado2i}
\left\{
\begin{array}{lcl}
\displaystyle
\frac{\partial \tilde{H}}{\partial t} = \xi_H \Delta \tilde{H} - \big[ \tau_H + l \gamma_H \Theta(\phi_1, \varphi_1)\big] \tilde{H} + \nu \tilde{A} - 
\vspace{0.2cm}
\\
\displaystyle
l \gamma_H \big[\Theta(\phi_1, \varphi_1) - \Theta(\phi_2, \varphi_2)\big]H_2  , &\textup{in}& Q,
\vspace{0.2cm}
\\
\displaystyle \frac{\partial \tilde{H}}{\partial \eta} (\cdot) =0, &\textup{on}& \Gamma,
\vspace{0.2cm}
\\
\displaystyle \tilde{H}(\cdot,0) = \tilde{H}_0(\cdot) = 0, &\textup{in}& \Omega.
\end{array}
\right.
\end{equation}

We observe that the system (\ref{separado1i}) satisfies the hypothesis of Proposition \ref{sol. Neumann}, therefore using the Lemmas \ref{PropertiesEtcFirst} and \ref{nãonegativa3}, $\tilde{A} \in W^{2,1}_4(Q)$ satisfies the following estimates
\begin{eqnarray*}
||\tilde{A}||_{W_{4}^{2, 1}(Q)} \leq C_p || l\beta_3\big[\Theta(\phi_1, \varphi_1) - \Theta(\phi_2, \varphi_2)\big] A_2||_{L^4(Q)} \\
\leq C_p \beta_3 k_A T^{\frac{1}{4}} |\Omega|^{\frac{1}{4}} C_1 ||(\phi_1, \varphi_1) - (\phi_2, \varphi_2)||_{L^{\infty}(Q)}.
\end{eqnarray*}

Thus, by Proposition \ref{icontLp01}, we obtain
\begin{equation}\label{ai}
  ||\tilde{A}||_{L^\infty(Q)} \leq C ||(\phi_1, \varphi_1) - (\phi_2, \varphi_2)||_{L^{\infty}(Q)},
\end{equation}
where $C$ depends on $C_p, \beta_3, k_A, T, \Omega, C_1$ and the immersion constant.

We observe that the system (\ref{separado2i}) satisfies the hypothesis of Proposition \ref{sol. Neumann}, therefore using the Lemma \ref{PropertiesEtcFirst}, $\tilde{H} \in W^{2,1}_4(Q)$ satisfies the following estimates
\begin{eqnarray*}
  ||\tilde{H}||_{W_{4}^{2, 1}(Q)} \leq C_p || \nu \tilde{A} - l \gamma_H \big[\Theta(\phi_1, \varphi_1) - \Theta(\phi_2, \varphi_2)\big]H_2||_{L^4(Q)} \\
  \leq C_p \nu ||\tilde{A}||_{L^4(Q)} + C_p \gamma_H T^{\frac{1}{4}} |\Omega|^{\frac{1}{4}} ||H_2||_{L^{\infty}(Q)} ||\Theta(\phi_1, \varphi_1) - \Theta(\phi_2, \varphi_2)||_{L^{\infty}(Q)} \\
  \leq C_p \nu ||\tilde{A}||_{L^\infty(Q)} +  C_p \gamma_H T^{\frac{1}{4}} |\Omega|^{\frac{1}{4}} C_1 ||H_2||_{L^{\infty}(Q)} ||(\phi_1, \varphi_1) - (\phi_2, \varphi_2)||_{L^{\infty}(Q)},
\end{eqnarray*}
therefore, using the estimate (\ref{ai}) and Proposition \ref{icontLp01}, we obtain
\begin{equation}\label{bx}
||\tilde{H}||_{L^\infty(Q)} \leq C ||(\phi_1, \varphi_1) - (\phi_2, \varphi_2)||_{L^{\infty}(Q)},
\end{equation}
where $C$ depends on $C_p$, $\nu$, $\gamma_H$, $T$, $\Omega$, $C_1$, $||H_2||_{L^\infty(Q)}$ and immersion constants.

Adding the inequalities obtained in (\ref{ai}) and (\ref{bx}), finally we get
\begin{eqnarray*}
  ||\tilde{A}||_{L^\infty(Q)} + ||\tilde{H}||_{L^\infty(Q)} &\leq& C ||(\phi_1, \varphi_1) - (\phi_2, \varphi_2)||_{L^{\infty}(Q)}.
\end{eqnarray*}

To show that $\Psi(l, \cdot, \cdot) $ is compact, we use the fact that immersion
$W^{2,1}_4(Q) \hookrightarrow L^\infty(Q)$ is compact and that
$\Psi(l, \cdot, \cdot)$ is the composition between the inclusion operator and the
solution operator, i.e., $\Psi(l, \cdot, \cdot): L^\infty(Q) \times L^\infty(Q)
\rightarrow W^{2,1}_4(Q) \times W^{2,1}_4(Q) \rightarrow L^\infty(Q) \times L^\infty(Q)$.
\hfill$\Box$
\\

\begin{lemma}
\label{unifcon}
Given a bounded subset $B \subset L^{\infty}(Q)\times L^{\infty}(Q)$, for each $(\phi, \varphi) \in B$, the mapping
$\Psi(\cdot, \phi, \varphi): [0, 1] \rightarrow L^{\infty}(Q)\times L^{\infty}(Q)$  is continuous, uniformly with respect to $B$.
\end{lemma}

\noindent {\bf Proof:}
Since $B \in L^\infty(Q) \times L^\infty(Q)$ is bounded, there is $r_B \ge 0$ such that, for any $(\phi, \varphi) \in B$ we have $||\phi||_{L^\infty(Q)} + ||\varphi||_{L^\infty(Q)} \leq r_B$. Now, let us fix $(\phi, \varphi) \in L^\infty(Q)$ and consider $l_1, l_2 \in [0,1]$ and denote $\Psi(\lambda_1, \phi, \varphi) = (A_1, H_1)$, $\Psi(l_2, \phi, \varphi) = (A_2, H_2)$; if $\tilde{A} = A_1 - A_2$ and $\tilde{H} = H_1 - H_2$, we have that $\tilde{A}$ and $\tilde{H}$ satisfy the following problems, respectively:
\begin{equation}
\label{B1}
\left\{
\begin{array}{lcl}
\displaystyle
\frac{\partial \tilde{A}}{\partial t} = \xi_A \Delta \tilde{A} + \bigg(r_A - (\mu_A + \epsilon_A) - \frac{r_A}{k_A}l_1(A_1 + A_2) - \beta_3\Theta(\phi, \varphi)l_1  \bigg)\tilde{A} -
\vspace{0.2cm}
\\
\displaystyle
\frac{r_A}{k_A}A_2^2(l_1 - l_2) - \beta_3\Theta(\phi, \varphi)A_2(l_1 - l_2), &\textup{in}& Q,
\vspace{0.2cm}
\\
\displaystyle \frac{\partial \tilde{A}}{\partial \eta} = 0, &\textup{on}& \Gamma,
\vspace{0.2cm}
\\
\displaystyle \tilde{A}(\cdot, 0) = \tilde{A}_0(\cdot) = 0, &\textup{in}& \Omega.
\end{array}
\right.
\end{equation}

\begin{equation}
\label{B2}
\left\{
\begin{array}{lcl}
\displaystyle
\frac{\partial \tilde{H}}{\partial t} = \xi_H \Delta \tilde{H} - (\tau_H + \gamma_H\Theta(\phi, \varphi)l_1)\tilde{H} +
\vspace{0.2cm}
\\
\nu\tilde{A} - \gamma_H\Theta(\phi, \varphi)H_2(l_1 - l_2), &\textup{in}& Q,
\vspace{0.2cm}
\\
\displaystyle \frac{\partial \tilde{H}}{\partial \eta} = 0, &\textup{on}& \Gamma,
\vspace{0.2cm}
\\
\displaystyle \tilde{H}(\cdot, 0) = \tilde{H}_0(\cdot) = 0, &\textup{in}& \Omega.
\end{array}
\right.
\end{equation}

We observe that the problem (\ref{B1}) satisfies the hypothesis of Proposition \ref{sol. Neumann}, therefore using the Lemmas \ref{PropertiesEtcFirst} and \ref{nãonegativa3}, $A \in W^{2,1}_4(Q)$ satisfies the following estimates
\begin{eqnarray*}
  ||\tilde{A}||_{W_{4}^{2, 1}(Q)} \leq C_p \bigg|\bigg|\frac{r_A}{k_A}A_2^2(l_1 - l_2) - \beta_3\Theta(\phi, \varphi)A_2(l_1 - l_2)\bigg|\bigg|_{L^4(Q)} \\
\leq  C_p k_A T^{\frac{1}{4}} |\Omega|^{\frac{1}{4}} (r_A +  \beta_3 ||N_0||_{L^\infty(\Omega)} + r_N T) |l_1 - l_2|.
\end{eqnarray*}

Thus, by Proposition \ref{icontLp01}, we obtain
\begin{equation}\label{aj}
  ||\tilde{A}||_{L^\infty(Q)} \leq C |l_1 - l_2|,
\end{equation}
where $C$ depends on $C_p$, $r_A$, $k_A$, $r_N$, $T$, $\Omega$, $\beta_3$, $||N_0||_{L^\infty(\Omega)}$ and the immersion constant.

We observe that the system (\ref{B2}) satisfies the hypothesis of Proposition \ref{sol. Neumann}, therefore using the Lemma \ref{PropertiesEtcFirst}, $\tilde{H} \in W^{2,1}_4(Q)$ satisfies the following estimates
\begin{eqnarray*}
||\tilde{H}||_{W_{4}^{2, 1}(Q)} \leq C_p ||\nu\tilde{A} - \gamma_H\Theta(\phi, \varphi)H_2(l_1 - l_2)||_{L^4(Q)} \\
\leq C_p \nu ||\tilde{A}||_{L^4(Q)} + C_p \gamma_H T^{\frac{1}{4}} |\Omega|^{\frac{1}{4}}||H_2||_{L^{\infty}(Q)} (||N_0||_{L^\infty(\Omega)} + r_N T) |l_1 - l_2| \\
\leq \bar{C}_p \nu ||\tilde{A}||_{L^\infty(Q)} +  C_p \gamma_H T^{\frac{1}{4}} |\Omega|^{\frac{1}{4}} ||H_2||_{L^{\infty}(Q)} (||N_0||_{L^\infty(\Omega)} + r_N T) |l_1 - l_2|,
\end{eqnarray*}
therefore, using the estimate (\ref{aj}) and Proposition \ref{icontLp01}, we obtain
\begin{equation}\label{Bii}
||\tilde{H}||_{L^\infty(Q)} \leq C |l_1 - l_2|.
\end{equation}
where $C$ depends on $C_p$, $r_A$, $k_A$, $r_N$, $\nu$, $\gamma_H$, $T$, $\Omega$, $\beta_3$, $||H_2||_{L^\infty(Q)}$, $||N_0||_{L^\infty(\Omega)}$ and immersion constants.

Adding the inequalities obtained in (\ref{aj}) and (\ref{Bii}), finally we get
\begin{eqnarray*}
  ||\tilde{A}||_{L^\infty(Q)} + ||\tilde{H}||_{L^\infty(Q)} &\leq& C |l_1 - l_2|.
\end{eqnarray*} 
\hfill$\Box$
\\

\begin{lemma}
\label{estimativa}
There exists a number $\rho > 0$ such that, for any $l \in [0,1]$ and any possible
fixed point $(A,H) \in  L^\infty(Q) \times L^\infty(Q)$ of $\Psi(l, \cdot, \cdot)$, there holds
$\|(A,H)\|_{L^\infty(Q)} < \rho$.
\end{lemma}

\noindent{\bf Proof:}
Let $(A,H) \in L^\infty(Q) \times L^\infty(Q)$ be a possible fixed point of $\Psi(l, \cdot, \cdot)$ associated to a given $l \in [0,1]$; then $(A,H)$ satisfies
\begin{equation}
\label{conex12}
\left\{
\begin{array}{lcl}
\displaystyle
\frac{\partial A}{\partial t} = \xi_A \Delta A + r_A A\left(1-\frac{A}{k_A}\right)-(\mu_A+\epsilon_A)A - l\beta_3 \Theta(\phi, \varphi) A  , &\textup{in}& Q,
\vspace{0.2cm}
\\
\displaystyle
\frac{\partial H}{\partial t} = \xi_H \Delta H + \nu A - \tau_H H - l \gamma_H \Theta(\phi, \varphi) H , &\textup{in}& Q,
\vspace{0.2cm}
\\
\displaystyle \frac{\partial A}{\partial \eta} (\cdot) =\frac{\partial H}{\partial \eta} (\cdot) =0, &\textup{on}& \Gamma,
\vspace{0.2cm}
\\
\displaystyle A(\cdot,0) = A_{0}(\cdot), H(\cdot,0) = H_0(\cdot), &\textup{in}& \Omega.
\end{array}
\right.
\end{equation}

By Proposition \ref{existenciaH}, we know that $H$ satisfies the estimate (\ref{existenciaH}), that is,
\begin{equation*}\label{HH}
  ||H||_{W^{2,1}_4(Q)} \leq C_p \bigg( \nu k_A T^{\frac{1}{4}}|\Omega|^{\frac{1}{4}} + ||H_0||_{W^{\frac{3}{2}}_4(\Omega)}\bigg).
\end{equation*}

Then, by Lemma \ref{icontLp01}, we have
\begin{equation*}\label{X}
  ||H||_{L^\infty(Q)} \leq \bar{C}_p \bigg( \nu k_A T^{\frac{1}{4}}|\Omega|^{\frac{1}{4}} + ||H_0||_{W^{\frac{3}{2}}_4(\Omega)}\bigg).
\end{equation*}

Therefore, by Lemma \ref{estimativa1}, we have
\begin{equation*}\label{X}
  ||A||_{L^\infty} + ||H||_{L^\infty(Q)} < \rho_1 + \bar{C}_p \bigg( \nu k_A T^{\frac{1}{4}}|\Omega|^{\frac{1}{4}} + ||H_0||_{W^{\frac{3}{2}}_4(\Omega)}\bigg),
\end{equation*}
which suggests we take
\begin{eqnarray*}
  \rho = \rho_1 + \bar{C}_p \bigg( \nu k_A T^{\frac{1}{4}}|\Omega|^{\frac{1}{4}} + ||H_0||_{W^{\frac{3}{2}}_4(\Omega)}\bigg).
\end{eqnarray*}
\hfill$\Box$
\\

\begin{lemma}
\label{fix}
The mapping $\Psi(0, \cdot, \cdot): L^\infty(Q) \times  L^\infty(Q) \rightarrow L^{\infty}(Q)\times L^{\infty}(Q)$ has a unique fixed point.
\end{lemma}

\noindent {\bf Proof:}
Indeed, letting $l =0$ in (\ref{P66}), $(A,H)$ is a fixed point of $\Psi(0, \cdot, \cdot)$ if, and only if, $(A,H)$ is the unique solution to the problem
\begin{equation*}
\label{P6}
\left\{
\begin{array}{lcl}
\displaystyle
\frac{\partial A}{\partial t} = \xi_A \Delta A + r_A A\left(1-\frac{A}{k_A}\right)-(\mu_A+\epsilon_A) A, &\textup{in}& Q,
\vspace{0.2cm}
\\
\displaystyle
\frac{\partial H}{\partial t} = \xi_H \Delta H + \nu A - \tau_H H, &\textup{in}& Q,
\vspace{0.2cm}
\\
\displaystyle \frac{\partial A}{\partial \eta} (\cdot)=\frac{\partial H}{\partial \eta} (\cdot) =0, &\textup{on}& \Gamma,
\vspace{0.2cm}
\\
\displaystyle A(\cdot,0) = A_{0}(\cdot), H(\cdot,0) = H_0(\cdot), &\textup{in}& \Omega.
\end{array}
\right.
\end{equation*}

Analogous to what was done in Subsection \ref{sec32}, we guarantee the existence of a unique solution $(A,H) \in  L^{\infty}(Q) \times L^\infty(Q)$ of this last problem; therefore $\Psi(0, \cdot, \cdot)$ has a unique fixed point in $L^\infty(Q) \times L^\infty(Q)$.
\hfill$\Box$
\\

\begin{proposition}
\label{h1}
There is a nonnegative solution $(\hat{A}, \hat{H}) \in W^{2,1}_4(Q) \times W^{2,1}_4(Q)$ of problem (\ref{P3}).
\end{proposition}

\noindent {\bf Proof:}
From Lemmas \ref{bemdefinido}, \ref{con}, \ref{unifcon}, \ref{estimativa} and \ref{fix}, we conclude that the mapping $\Psi: [0,1] \times L^\infty(Q) \times L^\infty(Q) \rightarrow L^\infty(Q) \times L^\infty(Q)$ satisfies the hipotheses of the Leray-Schauder's fixed point theorem (see Friedman \cite[pp.~189, Theorem 3]{Friedman}). Thus, there exists $(\hat{A}, \hat{H}) \in L^\infty(Q) \times L^\infty(Q)$ such that $\Psi(1, \hat{A}, \hat{H}) = (\hat{A}, \hat{H})$. Moreover, by Lemmas \ref{nãonegativa1}, \ref{existenciaA}, \ref{existenciaH} and \ref{nãonegativa2}, $(\hat{A}, \hat{H}) \in W^{2,1}_4(Q) \times W^{2,1}_4(Q)$ is nonnegative, meets the estimates (\ref{AA}), (\ref{existenciaH}) and $(\hat{A}, \hat{H})$ is the required solution of (\ref{P3}).
\hfill$\Box$
\\

\section{Proof of Theorem \ref{Teorema1}}

\begin{proposition}
\label{h2}
There is a nonnegative solution $(\hat{N}, \hat{A}, \hat{H}) \in L^\infty(Q) \times W^{2,1}_4(Q) \times W^{2,1}_4(Q)$ of the modified problem (\ref{P01}).
\end{proposition}

\noindent {\bf Proof:}
Just combine the Proposition \ref{h1}, the Remark \ref{obs1} and the Lemma \ref{PropertiesEtcFirst}.
\hfill$\Box$
\\

\begin{remark}
\label{estimativaN}
We affirm that $\hat{N} \in W$. Indeed, of the Lemma \ref{PropertiesEtcFirst} we know that $\hat{N} = \Theta(\hat{H}, \hat{A}) \in L^\infty(Q)$. Moreover, returning to the first equation of (\ref{P01}), using the Lemmas \ref{PropertiesEtcFirst} and \ref{nãonegativa3} and the fact that $\hat{H} \in W^{2,1}_4(Q) \subset L^\infty(Q)$, it follow that:
\begin{equation}\label{Nt}
  \bigg|\frac{\partial \hat{N}}{\partial t}\bigg| \leq r_N + (\mu_N + \beta_1 k_A + \alpha_H\gamma_H ||\hat{H}||_{L^\infty(Q)})  (||N_0||_{L^\infty(\Omega)} + r_N T),
\end{equation}
a.e. in $Q$, i.e., $\hat{N}_t \in L^\infty(Q)$.
\end{remark}

\begin{proposition}
\label{j1}
There is a nonnegative solution $(N, A, H) \in W \times W^{2,1}_4(Q) \times W^{2,1}_4(Q)$ of problem (\ref{0riginalEquations}).
\end{proposition}

\noindent {\bf Proof:}
Just combine the Proposition \ref{h2} and the Remarks \ref{obs2} and \ref{estimativaN}.
\hfill$\Box$
\\

\begin{proposition}
\label{j2}
The solution $(N, A, H)$ of the problem (\ref{0riginalEquations}) is unique.
\end{proposition}

\noindent {\bf Proof:}
Let $(N_1, A_1, D_1)$ and $(N_2, A_2, H_2)$ be solutions to the problem (\ref{0riginalEquations}); if $\tilde{N} = N_1 - N_2, \tilde{A} = A_1 - A_2$ and $\tilde{H} = H_1 - H_2$, then $\tilde{N}$, $\tilde{A}$ and $\tilde{H}$ satisfy the following problems, respectively:
\begin{equation}
\label{original1}
\left\{
\begin{array}{lcl}
\displaystyle
\frac{\partial \tilde{N}}{\partial t} = - \mu_N \tilde{N} - \beta_1 N_1 \tilde{A} - \beta_1 A_2 \tilde{N} - \alpha_H\gamma_H N_1 \tilde{H} - \alpha_H\gamma_H H_2 \tilde{N} , & \textup{in}& Q,
\vspace{0.2cm}
\\
\displaystyle
\displaystyle \tilde{N}(\cdot,0) = \tilde{N}_0(\cdot) = 0, &\textup{in}& \Omega,
\end{array}
\right.
\end{equation}

\begin{equation}
\label{original2}
\left\{
\begin{array}{lcl}
\displaystyle
\frac{\partial \tilde{A}}{\partial t} = \xi_A \Delta \tilde{A} + r_A \tilde{A} - \frac{r_A}{k_A}(A_1 + A_2)\tilde{A} -(\mu_A+\epsilon_A)\tilde{A} - \beta_3 N_1 \tilde{A} - \beta_3 A_2\tilde{N}, &\textup{in}& Q,
\vspace{0.2cm}
\\
\displaystyle \frac{\partial \tilde{A}}{\partial \eta} (\cdot) =0, &\textup{on}& \Gamma,
\vspace{0.2cm}
\\
\displaystyle \tilde{A}(\cdot,0) = \tilde{A}_0(\cdot) = 0, &\textup{in}& \Omega,
\end{array}
\right.
\end{equation}

\begin{equation}
\label{original3}
\left\{
\begin{array}{lcl}
\displaystyle
\frac{\partial \tilde{H}}{\partial t} = \xi_H \Delta \tilde{H} + \nu \tilde{A} - \tau_H \tilde{H} - \gamma_H N_1 \tilde{H} - \gamma_H H_2 \tilde{N} , &\textup{in}& Q,
\vspace{0.2cm}
\\
\displaystyle \frac{\partial \tilde{H}}{\partial \eta} (\cdot) =0, &\textup{on}& \Gamma,
\vspace{0.2cm}
\\
\displaystyle \tilde{H}(\cdot,0) = \tilde{H}_0(\cdot) = 0, &\textup{in}& \Omega.
\end{array}
\right.
\end{equation}

Multiplying the first equation of (\ref{original1}) by $\tilde{N}$, integrating into $\Omega$, using the fact that $N_1 \leq ||N_0||_{L^\infty(\Omega)} + r_N T$ and the inequality of Young, we have
\begin{eqnarray*}
  \frac{1}{2} \frac{d}{dt} \int_{\Omega} \tilde{N}^2 dx &=& - \mu_N \int_{\Omega} \tilde{N}^2 dx - \beta_1 \int_{\Omega} N_1 \tilde{A} \tilde{N} dx - \beta_1 \int_{\Omega} A_2 \tilde{N}^2 dx \\
  &-& \alpha_H\gamma_H \int_{\Omega} N_1 \tilde{H} \tilde{N} dx - \alpha_H\gamma_H \int_{\Omega} H_2 \tilde{N}^2 dx \\
  &\leq& (||N_0||_{L^\infty(\Omega)})\bigg(\beta_1 \int_{\Omega} |\tilde{A}||\tilde{N}| dx + \alpha_H\gamma_H \int_{\Omega} |\tilde{H}||\tilde{N}| dx\bigg) \\
  &\leq& C \int_{\Omega} (\tilde{A}^2 + \tilde{N}^2 + \tilde{H}^2) dx,
\end{eqnarray*}
where $C$ depends on $\beta_1$, $\alpha_H$, $\gamma_H$, $r_N$, $T$ and $||N_0||_{L^\infty(\Omega)}$.

Now, multiplying the first equation of (\ref{original2}) by $\tilde{A}$, integrating into $\Omega$, using the fact that $A \leq k_A$ and the inequality of Young, we obtain
\begin{eqnarray*}
  \frac{1}{2} \frac{d}{dt} \int_{\Omega} \tilde{A}^2 dx &=& - \xi_A \int_{\Omega} |\nabla \tilde{A}|^2 dx + r_A \int_{\Omega} \tilde{A}^2 dx - \frac{r_A}{k_A} \int_{\Omega}(A_1 + A_2)\tilde{A}^2 dx  \\
  &-& (\mu_A+\epsilon_A)\int_{\Omega}\tilde{A}^2 dx - \beta_3 \int_{\Omega} N_1 \tilde{A}^2 dx - \beta_3 \int_{\Omega} A_2\tilde{N} \tilde{A} dx \\
  &\leq& r_A \int_{\Omega} \tilde{A}^2 dx + \beta_3 k_A \int_{\Omega} |\tilde{N}| |\tilde{A}| dx \\
  &\leq& C  \int_{\Omega} (\tilde{A}^2 + \tilde{N}^2 + \tilde{H}^2) dx,
\end{eqnarray*}
where $C$ depends on $r_A$, $k_A$ and $\beta_3$.

Lastly, multiplying the first equation of (\ref{original3}) by $\tilde{H}$, integrating into $\Omega$, using the fact that $A \leq k_A$ and the inequality of Young, we obtain
\begin{eqnarray*}
  \frac{1}{2} \frac{d}{dt} \int_{\Omega} \tilde{H}^2 dx &=& - \xi_H \int_{\Omega} |\nabla \tilde{H}|^2 dx + \nu \int_{\Omega} \tilde{A} \tilde{H} dx - \tau_H \int_{\Omega} \tilde{H}^2 dx \\
  &-& \gamma_H \int_{\Omega} N_1 \tilde{H}^2 dx - \gamma_H \int_{\Omega} H_2 \tilde{N} \tilde{H} dx \\
  &\leq& \nu \int_{\Omega} |\tilde{A}| |\tilde{H}| dx + \gamma_H ||H_2||_{L^\infty(Q)} \int_{\Omega} |\tilde{N}| |\tilde{H}| dx \\
  &\leq& C \int_{\Omega} (\tilde{A}^2 + \tilde{N}^2 + \tilde{H}^2) dx,
\end{eqnarray*}
where $C$ depends on $\nu$, $\gamma_H$ and $||H_2||_{L^\infty(Q)}$.

Thus,
\begin{eqnarray*}
  \frac{d}{dt}\bigg( \int_{\Omega} (|\tilde{N}|^2 + |\tilde{A}|^2 + |\tilde{H}|^2) dx \bigg) &\leq& C \int_{\Omega} (|\tilde{N}|^2 + |\tilde{A}|^2+ |\tilde{D}|^2) dx
\end{eqnarray*}
and using the Gronwall's inequality, we obtain
\begin{equation*}
\int_{\Omega} (|\tilde{N}|^2 + |\tilde{A}|^2 + |\tilde{H}|^2) dx \leq e^{CT}  \int_{\Omega} (|\tilde{N}_0|^2 + |\tilde{A}_0|^2 + |\tilde{D}_0|^2) dx = 0,
\end{equation*}
that is, $||\tilde{N}(\cdot, t)||_{L^2(\Omega)}^2 + ||\tilde{A}(\cdot, t)||_{L^2(\Omega)}^2 + ||\tilde{D}(\cdot, t)||_{L^2(\Omega)}^2 = 0$, for all $t \in (0, T)$, where we conclude $\tilde{N} = \tilde{A} = \tilde{D} = 0$ a.e. in $Q$ and therefore $N_1 = N_2, A_1 = A_2$ and $D_1 = D_2$ a.e. in $Q$.
\hfill$\Box$
\\

\section{Numerical simulations}

In this section, we illustrate model behavior with numerical simulations. The simulations settings are the following. We consider the spatial domain as a square, $\Omega=[0,L] \times [0,L]$, with $L=1$, discritised with steps $\Delta x = \Delta y = 0.01$. The coupled ODEs arising from the spatial discretisation are solved with the method of lines in \textit{Mathematica}. The simulations run from time $t=0$ until $t=50$.

To avoid large numbers and numerical instabilities, we rescale the populations with respect to their possible maximum values, setting $N \leftarrow N/(r_N/\mu_N)$ and $A\leftarrow A/K_A$. Therefore, the population sizes range from $0$ to $1$. The parameter values used to simulate the model were
\[
r_N=1, \
\mu_N=1, \
\beta_3=1, \
r_A = 1, \
K_A = 1, \
\beta_1 = 1.5, \
\mu_A =0.05, \
\epsilon_A =0.05, \
\]
\[
\nu= 2, \
\tau_H=0.9, \
\gamma_H=0.01, \
\alpha_H =2200, \
\xi_H=0.01,
\xi_A=0.001.
\]
These values were chosen to describe: normal cells that reach the equilibrium $N=r_N/\mu_N=1$ at absence of tumor cells; a tumor with the same carrying capacity of normal cells ($K_A=r_N/\mu_N$), and that causes more damage to normal cells than the contrary ($\beta_1>\beta_3$) but is not able to invade the tissue without production of lactic acid (see simulation 1 below); a faster lactic acid diffusion in comparison with tumor cells ($\xi_H>\xi_A$); and a acid damage high enough to allow tumor invasion ($\alpha_H=2200$, see simulations below).

The initial conditions for numerical simulations were $N(x,0)=r_N/\mu_N$,
\[
A(x,0)=A_0 \exp \left(-\delta_A \left((x_1-L/2)^2+(x_2-L/2)^2\right)\right),
\]
and $H(x,0)=0$, with $x=(x_1,x_2)\in \Omega \subset \mathbb{R}^2$, $A_0=0.22$, $\delta_A=1000$. These initial condition describe the normal cells at the tissue normal homoeostatic state and the onset of a very small tumor in the middle of the tissue, with zero initial concentration of lactic acid (see Figure 1).

We present the following simulation results. In the first simulation (Figure 1), we set the lactic acid production to zero ($\nu=0$), to observe the extinction of tumor cells with no ability to produce lactic acid. In the second simulation, we set the parameters as stated above and observe the invasion of a tumor and substantial reduction of normal cells (Figure 2). In the third simulation, we replace the Neuman homogeneus conditions on tumor cells and lactic acid by Dirichlet homogeneous conditions, representing a situation where the tissue boundary ($\partial \Omega$) is hostile to the tumor (Figure 3). We observe a similar behavior, with exception that in the boundary the normal cells remain at high numbers, due to the low presence of tumor cells and lactic acid.

\begin{figure}[!htb]
\centering
\includegraphics[width=0.99\linewidth]{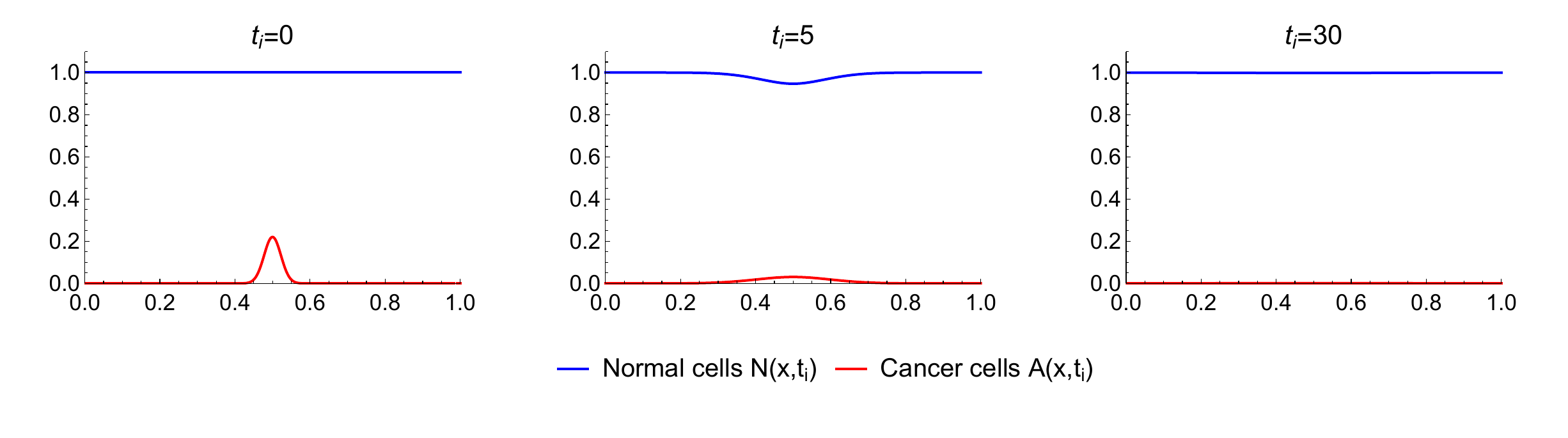}
\caption{Simulation results of model (\ref{0riginalEquations}) with no lactic acid production. The plots correspond to transversal sections of model solutions $A(x,t)$ (cancer cells, red) and $N(x,t)$ (normal cells, blue) at $x=(x_1,x_2)=(L/2,x_2)$, $x_2 \in [0,L]$, at time points $t=0,5,30$. The initial tumor cells are extinct, due the negative effect exerted by the normal cells.
}
\end{figure}

\begin{figure}[!htb]
\centering
\includegraphics[width=0.24\linewidth]{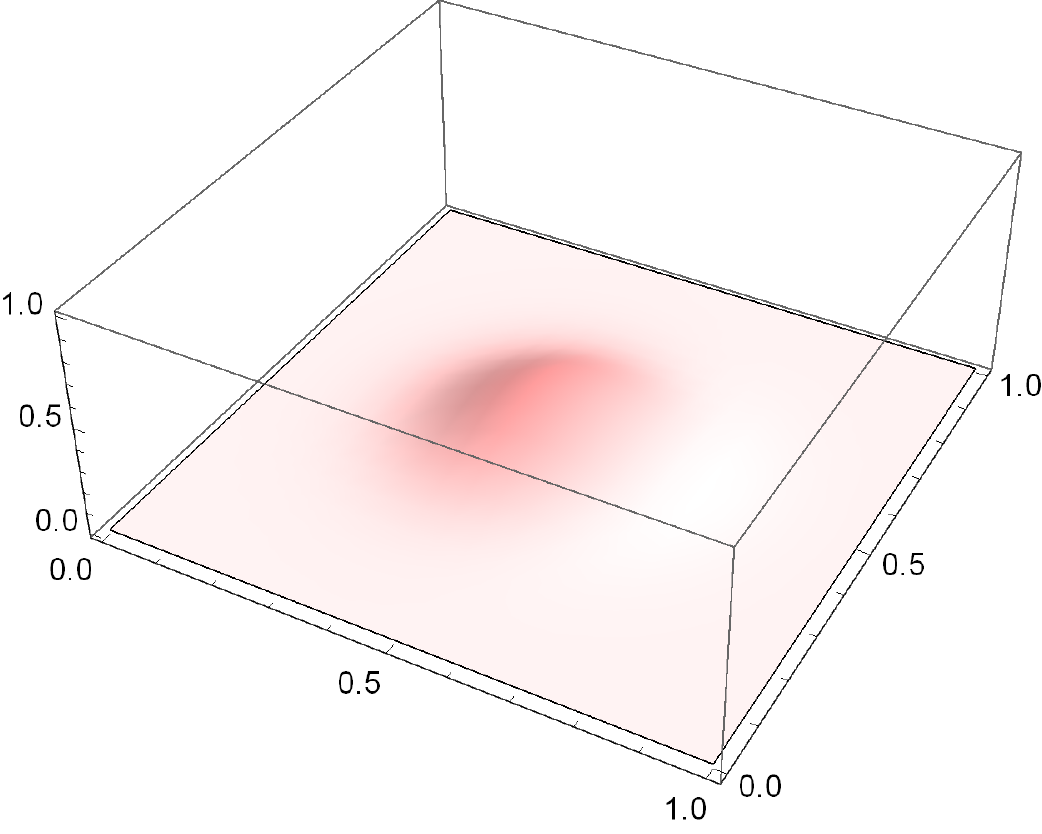}
\includegraphics[width=0.24\linewidth]{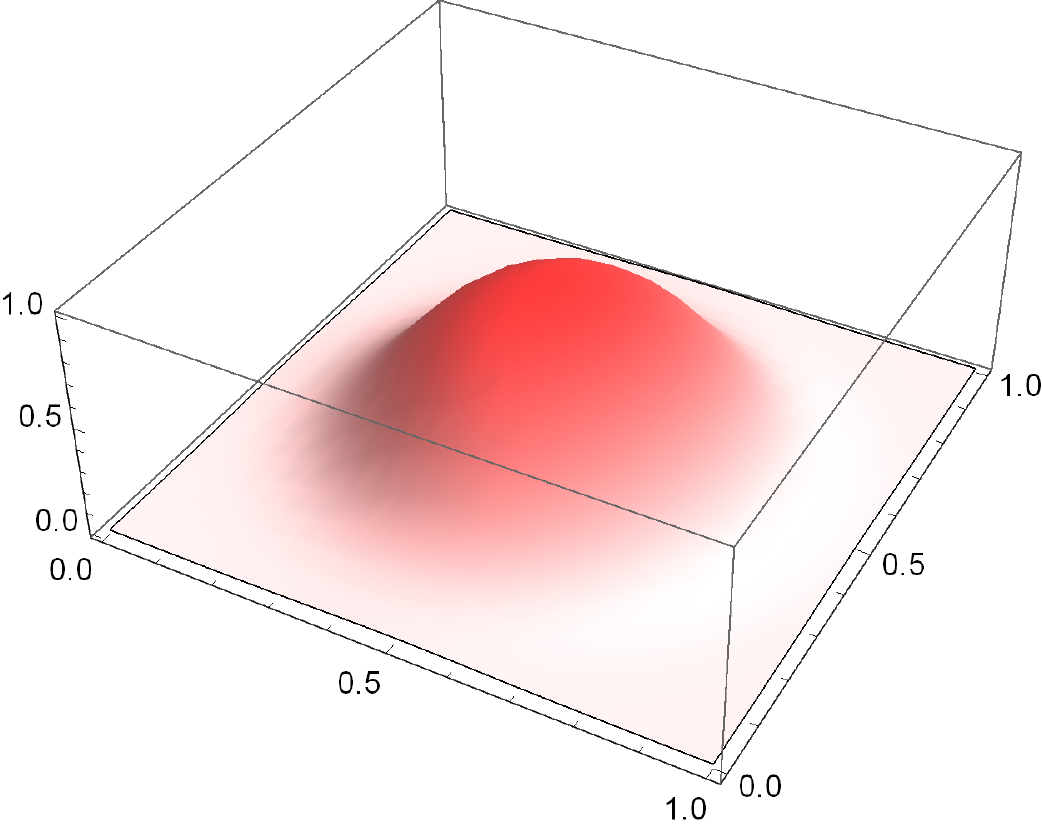}
\includegraphics[width=0.24\linewidth]{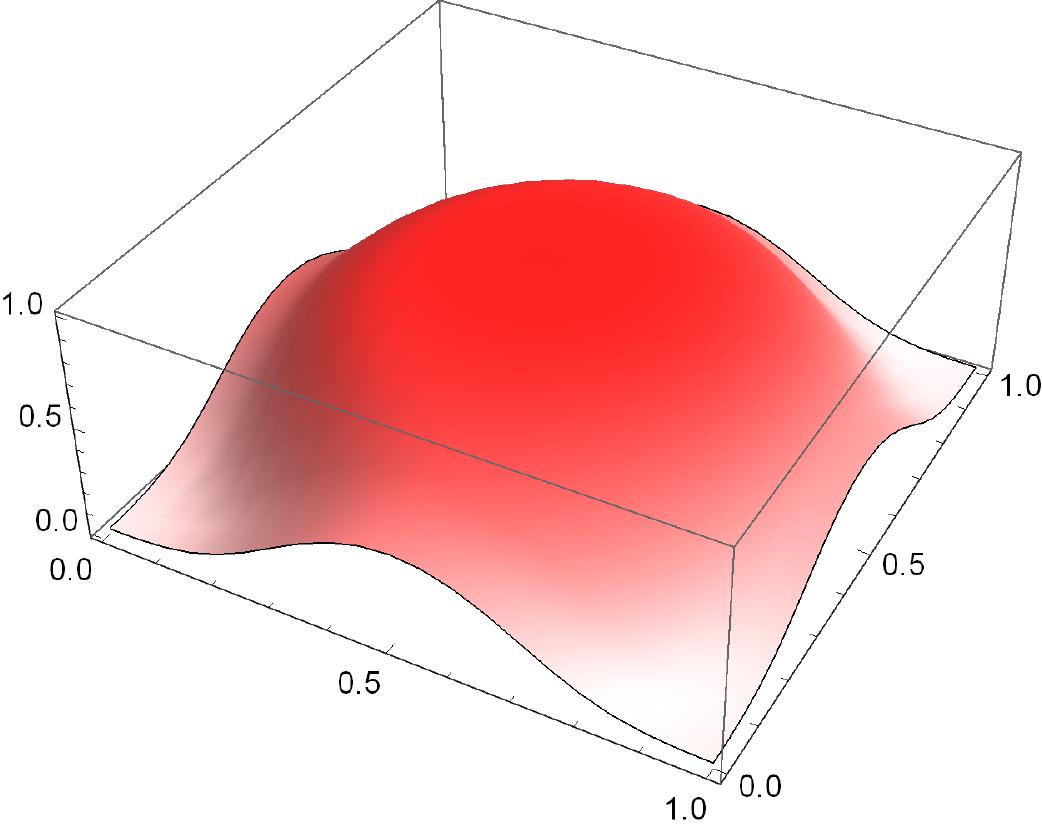}
\includegraphics[width=0.24\linewidth]{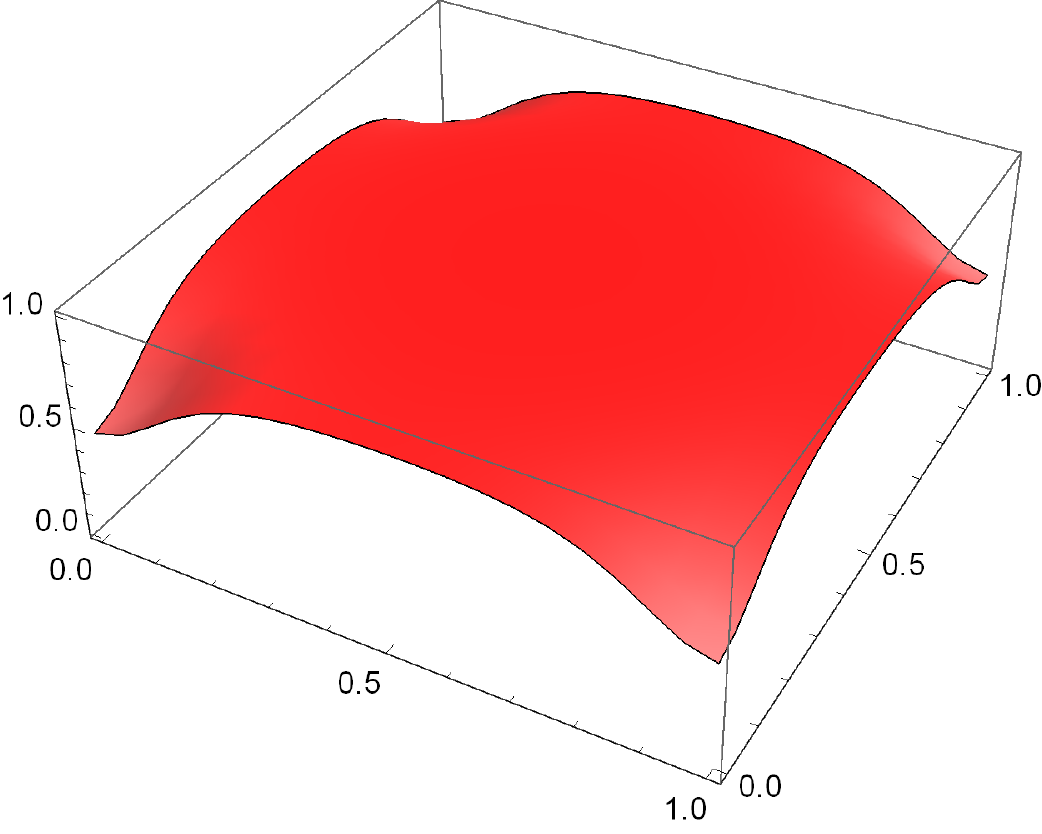}
\includegraphics[width=0.24\linewidth]{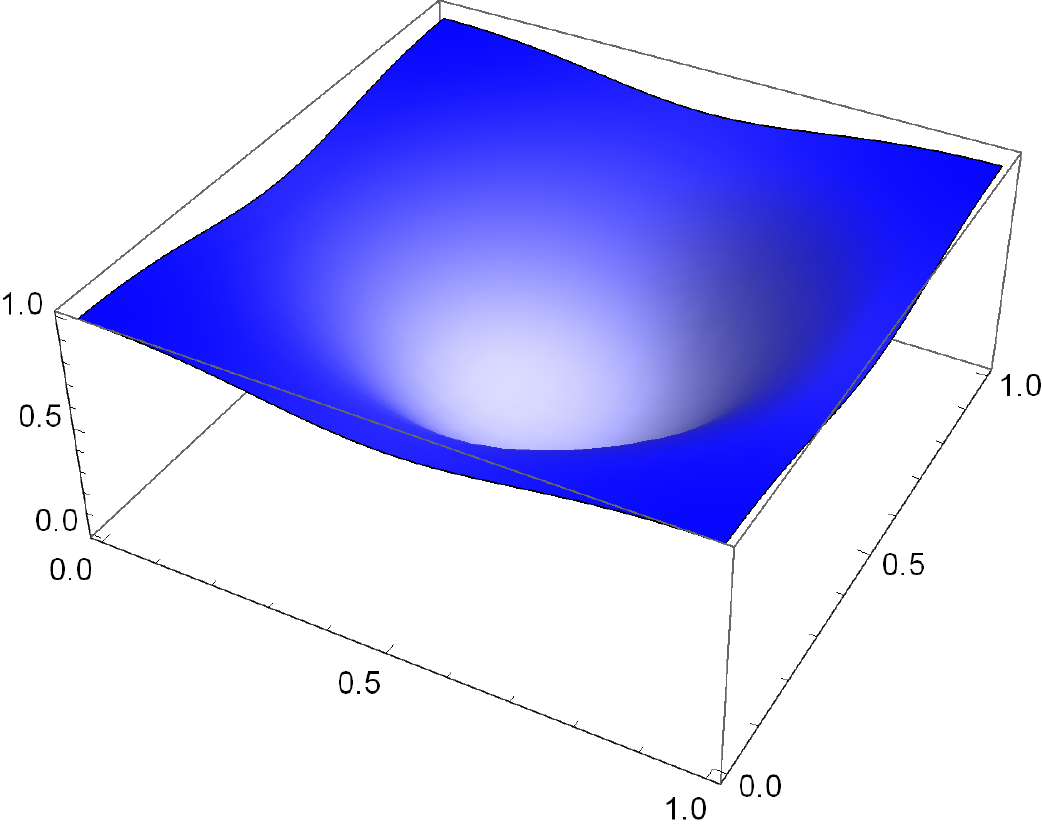}
\includegraphics[width=0.24\linewidth]{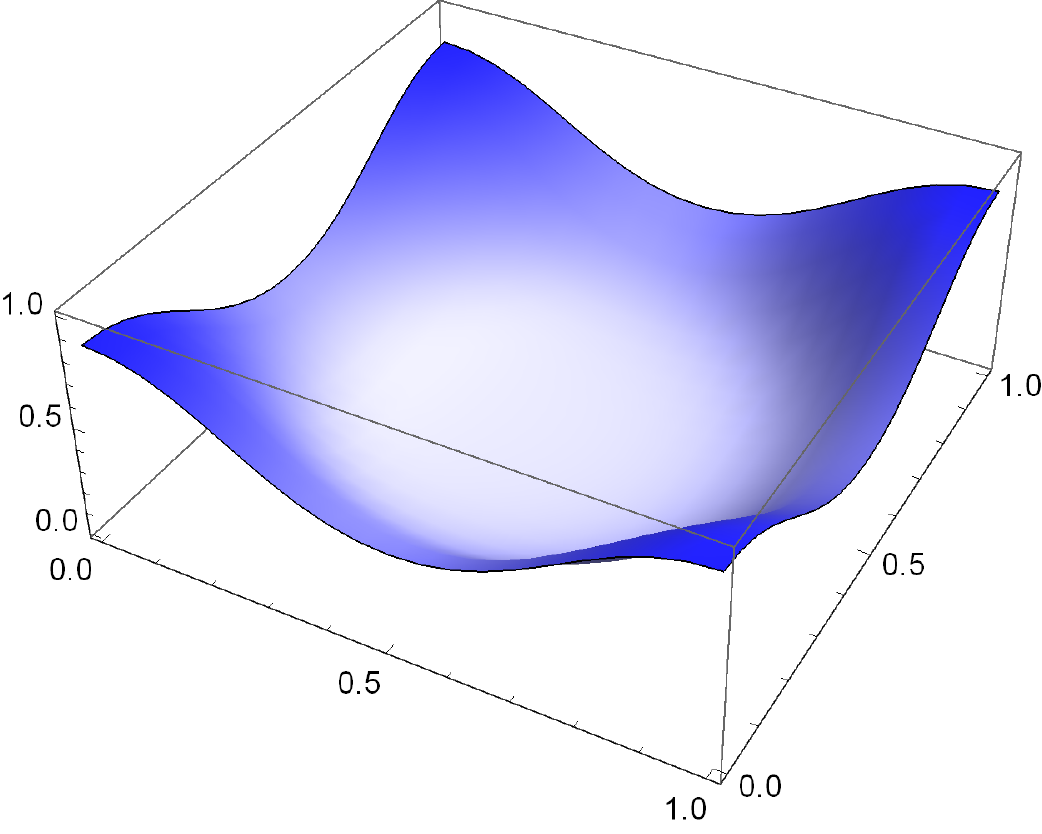}
\includegraphics[width=0.24\linewidth]{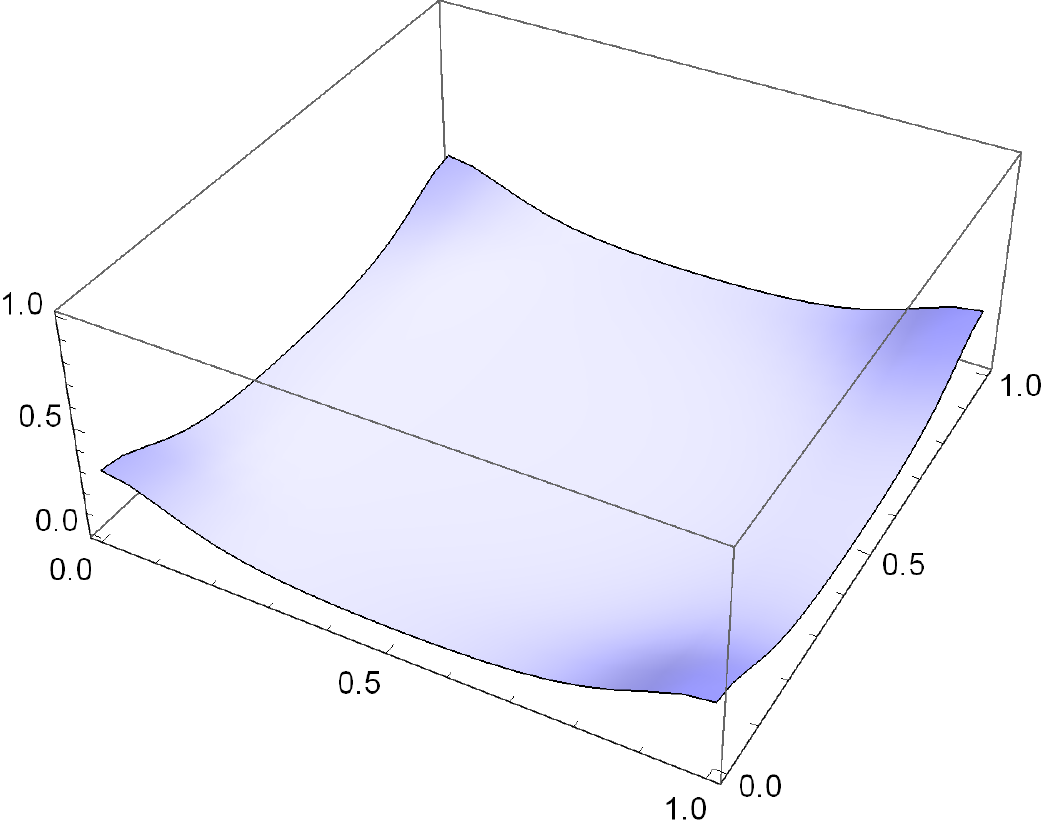}
\includegraphics[width=0.24\linewidth]{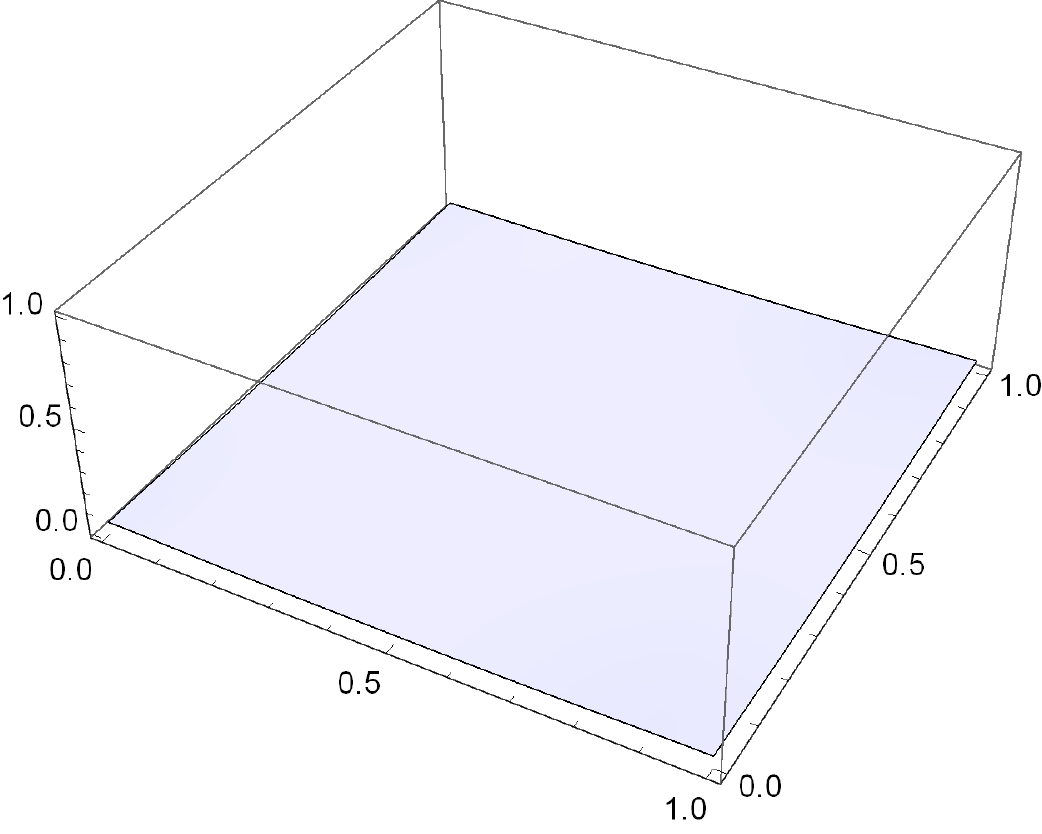}
\includegraphics[width=0.24\linewidth]{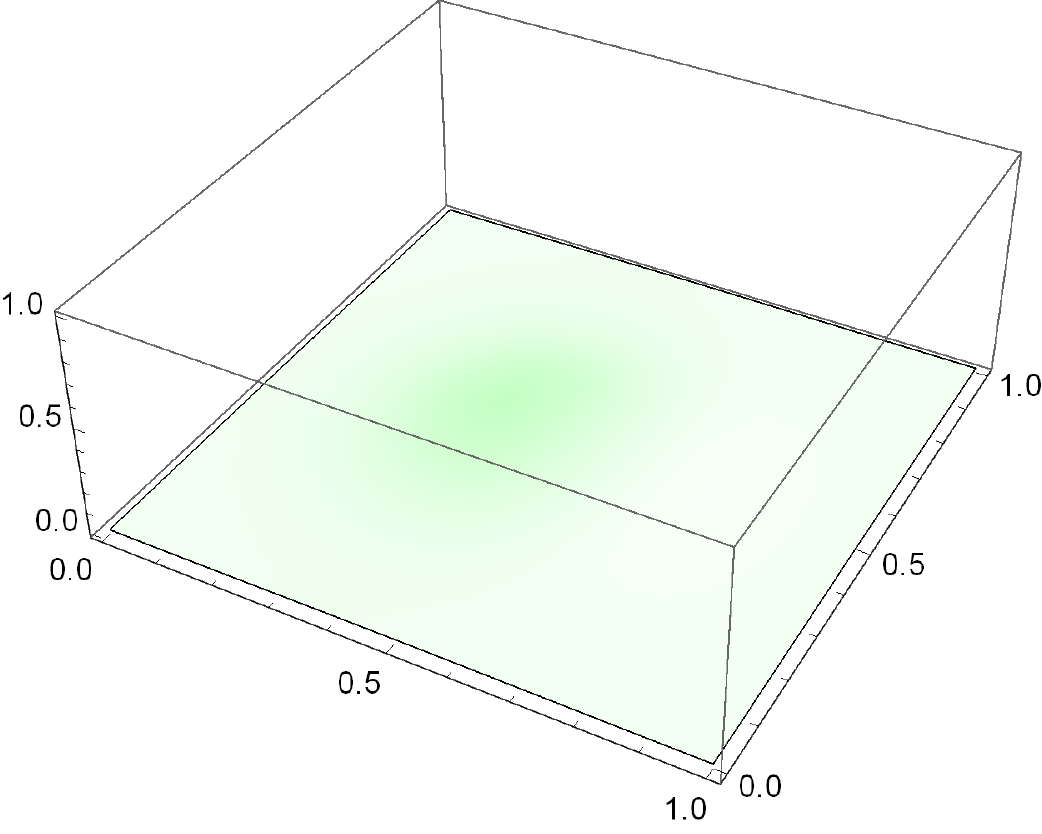}
\includegraphics[width=0.24\linewidth]{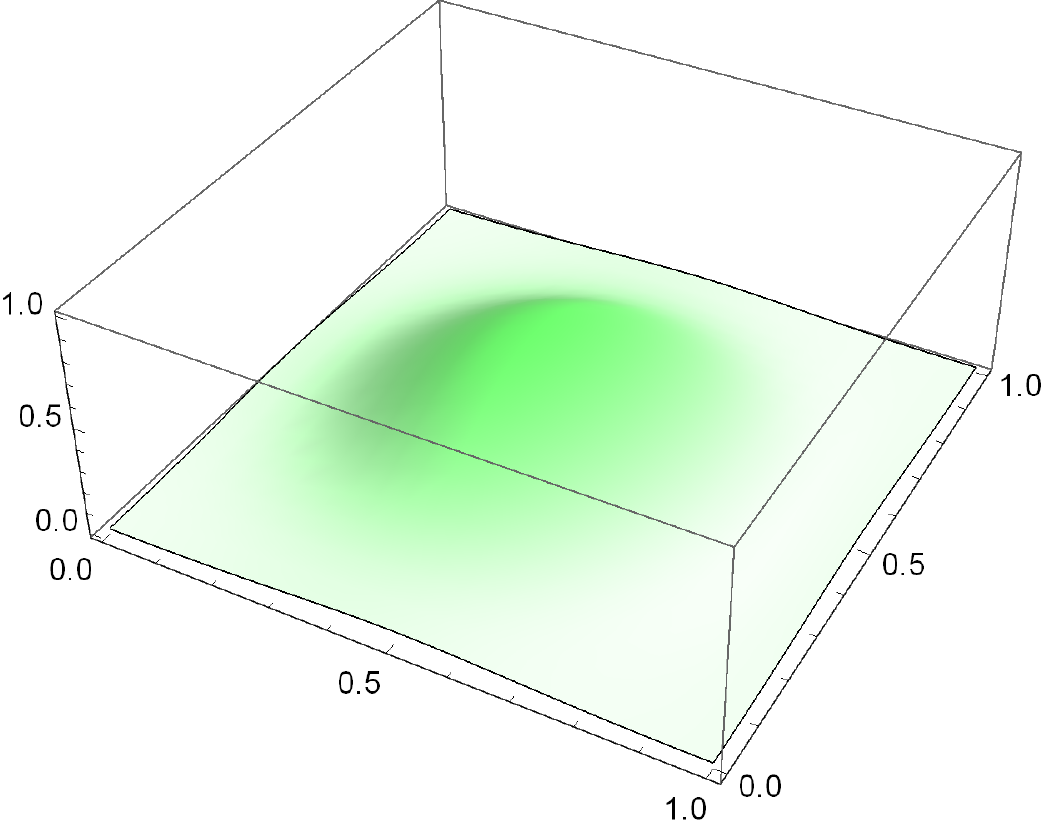}
\includegraphics[width=0.24\linewidth]{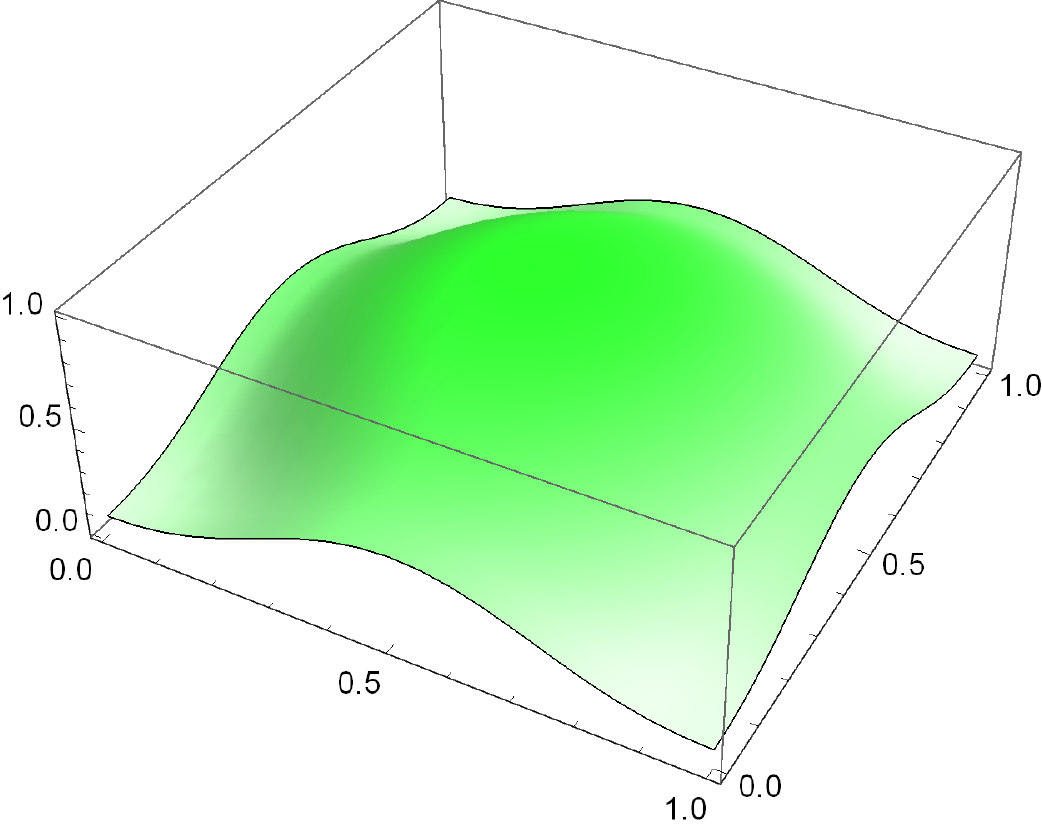}
\includegraphics[width=0.24\linewidth]{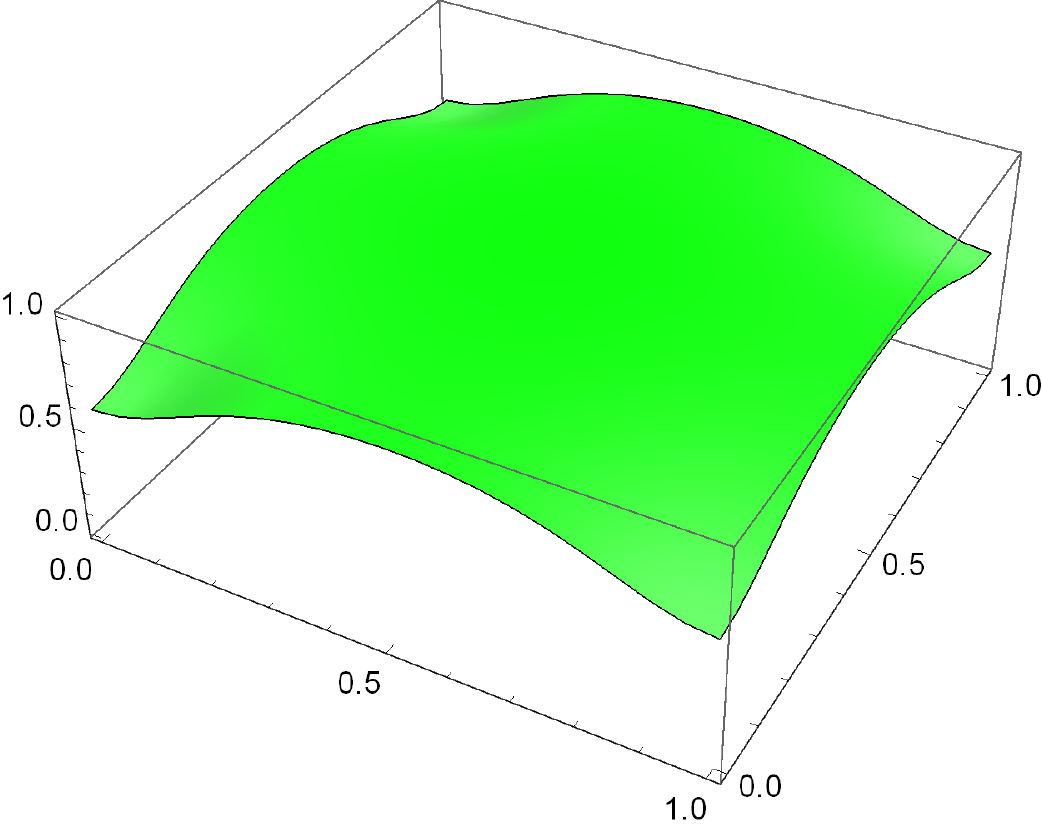}
\caption{Simulation results of model (\ref{0riginalEquations}). Plots of model solutions $A(x,t)$ (cancer cells, red, top row), $N(x,t)$ (normal cells, blue, middle row) and $H(x,t)$ (lactic acid concentration, green, bottom row) at time points $t=23, 27,31,35$, with $x\in \Omega=[0,L]\times [0,L]$. The initial tumor cells survive, produce lactic acid and invade the tissue, leading the drastic reduction on the number of normal cells.
}
\end{figure}

\begin{figure}[!htb]
\centering
\includegraphics[width=0.24\linewidth]{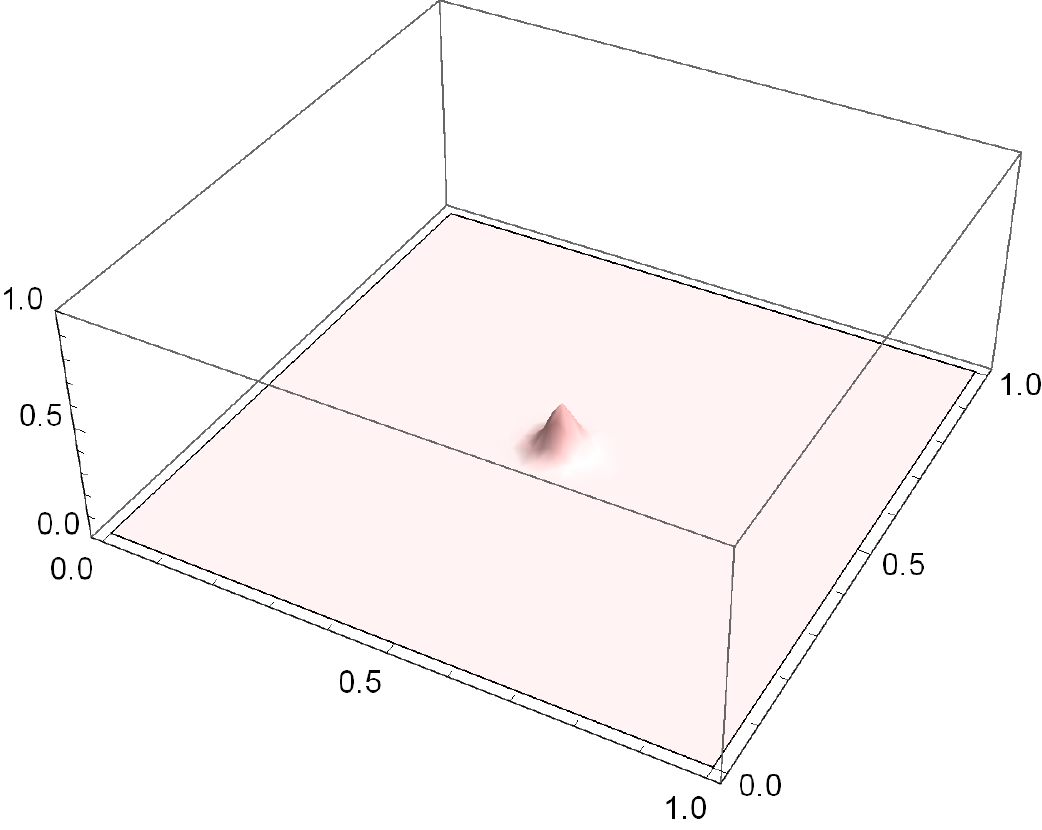}
\includegraphics[width=0.24\linewidth]{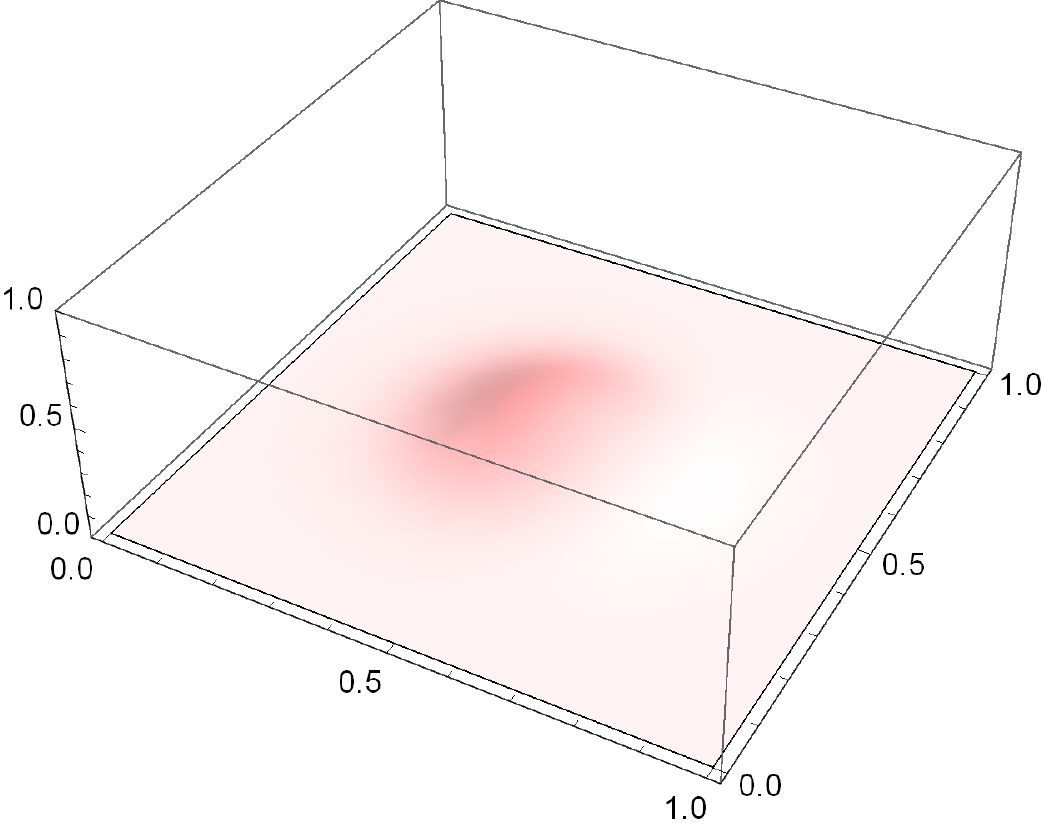}
\includegraphics[width=0.24\linewidth]{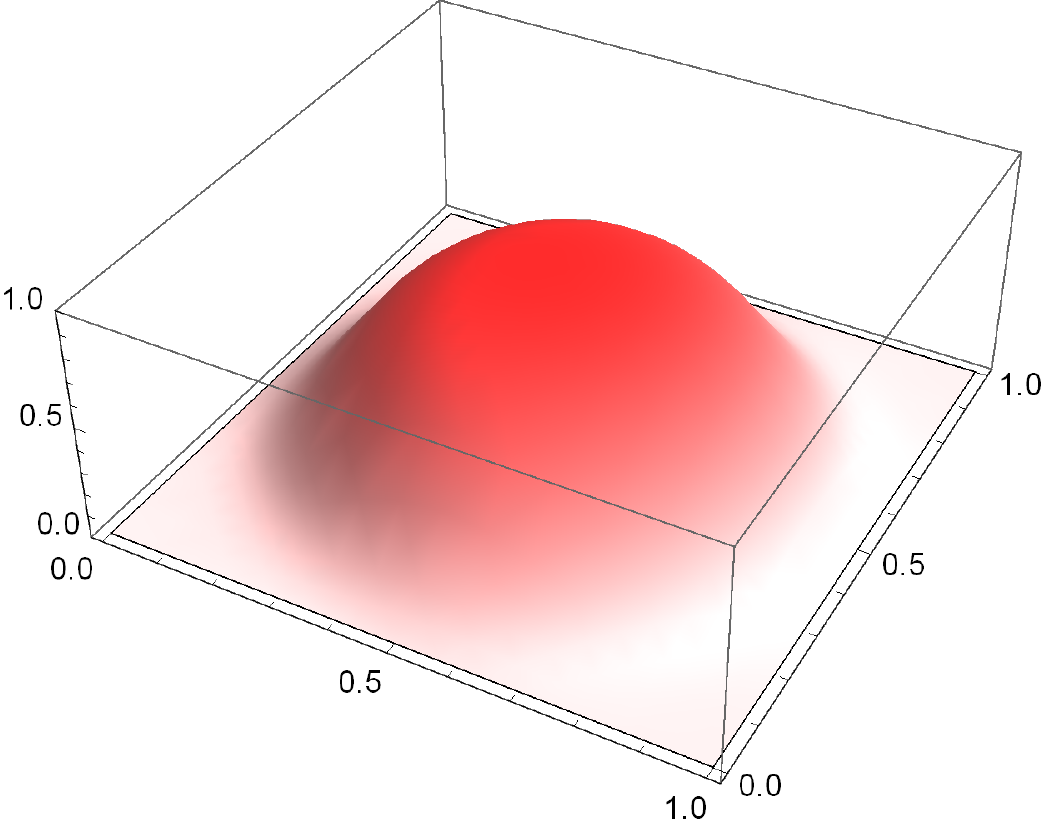}
\includegraphics[width=0.24\linewidth]{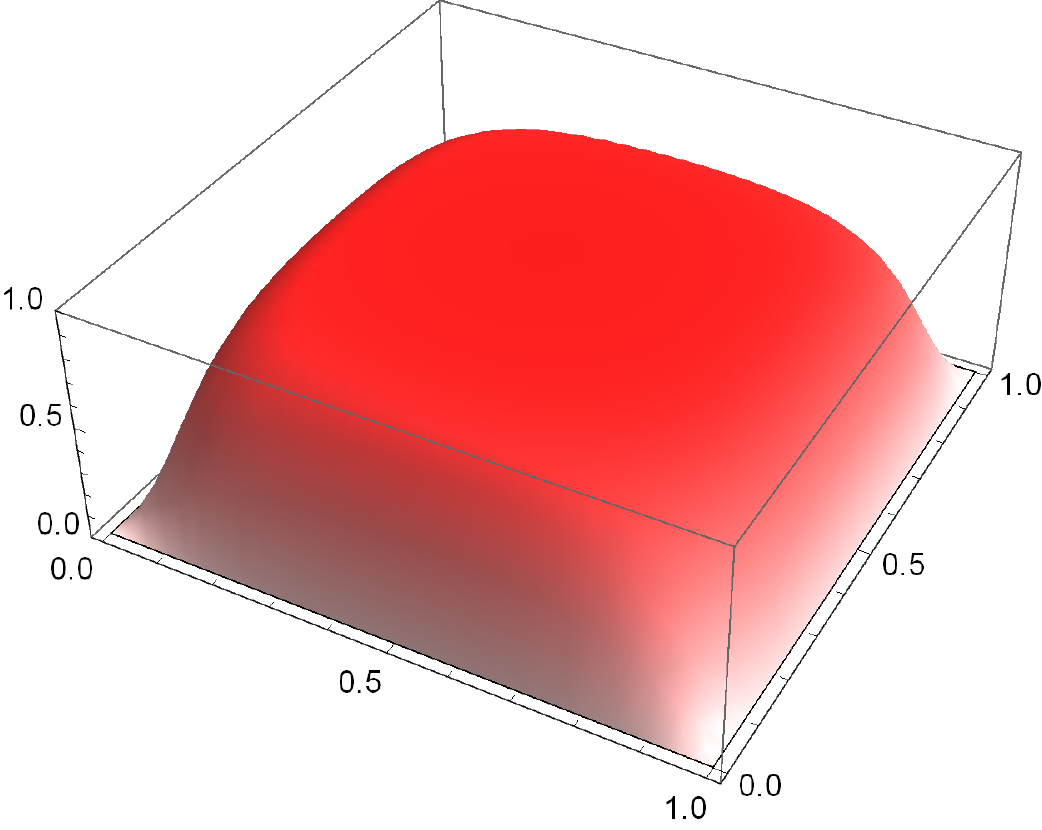}
\includegraphics[width=0.24\linewidth]{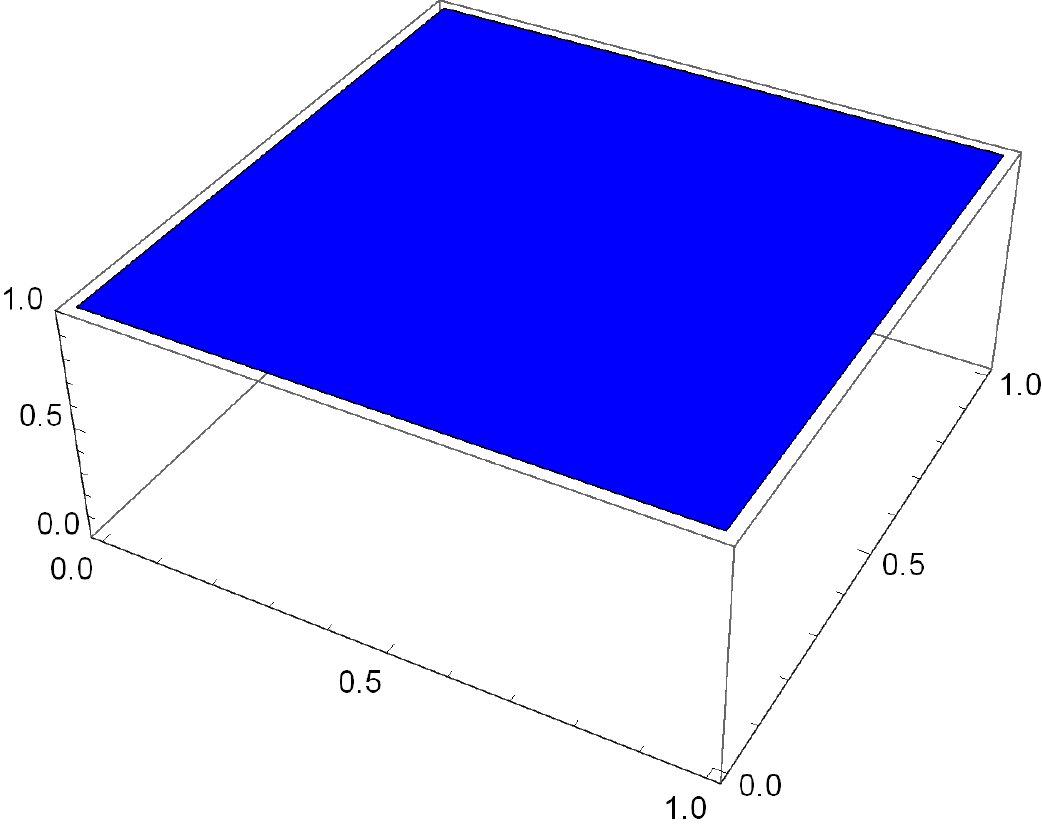}
\includegraphics[width=0.24\linewidth]{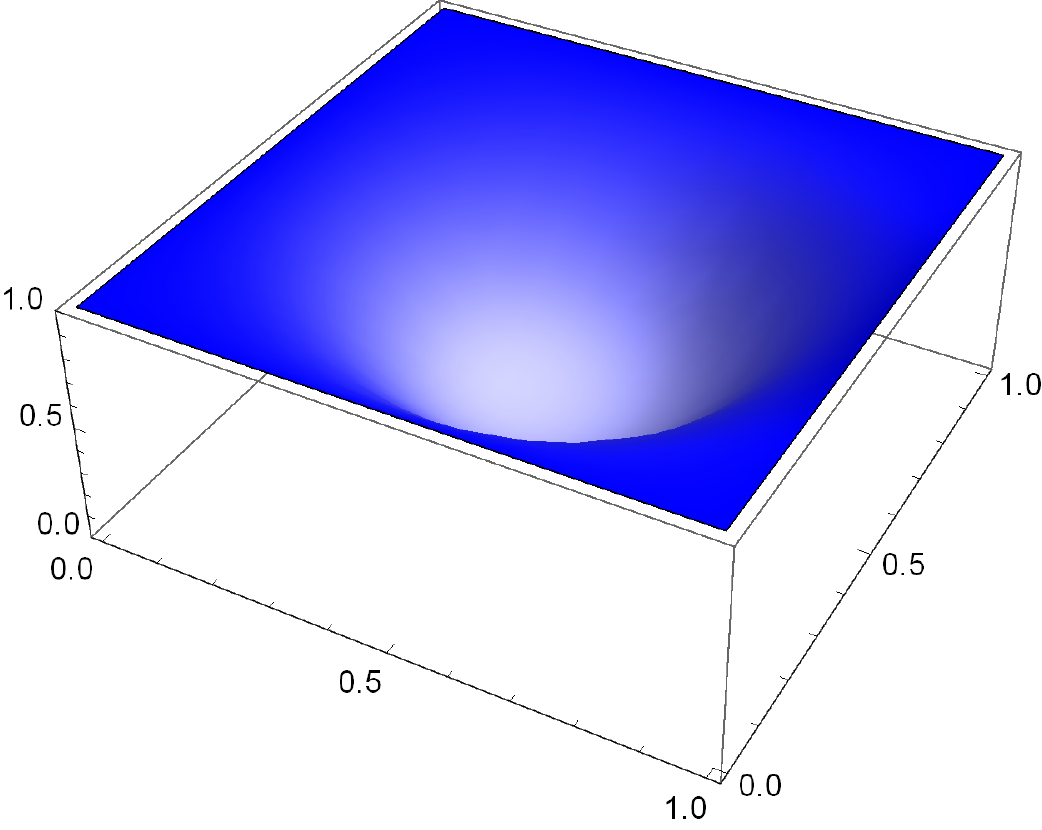}
\includegraphics[width=0.24\linewidth]{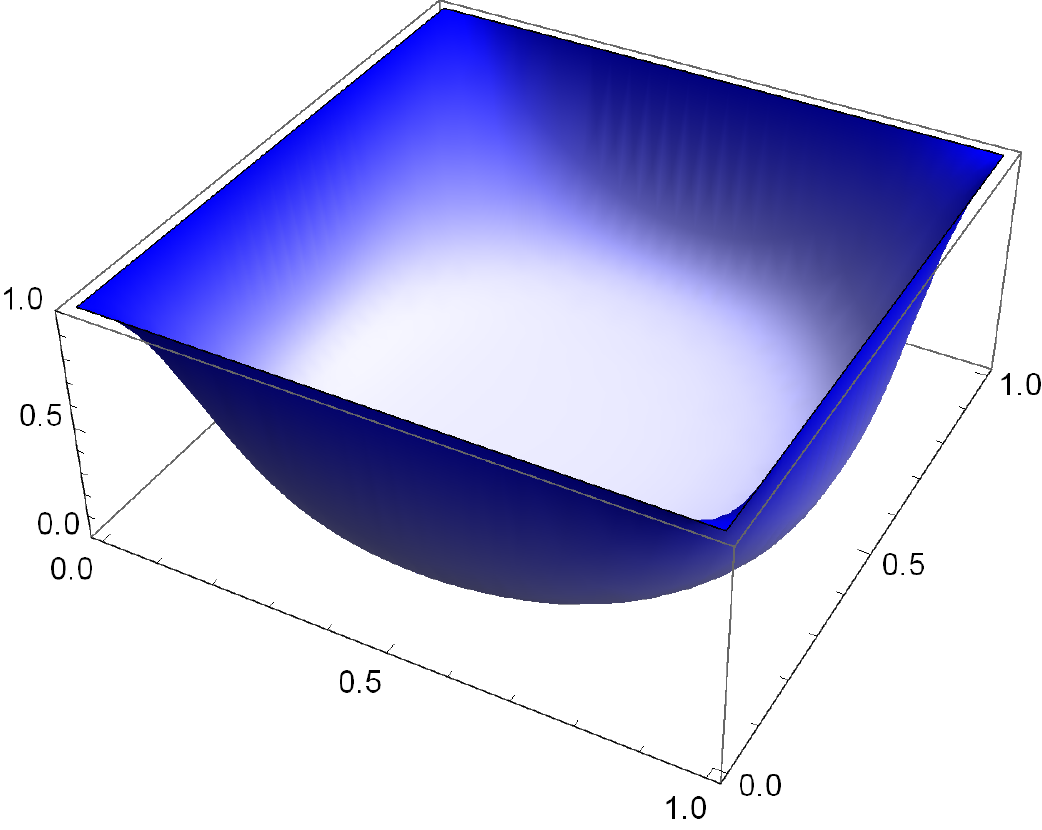}
\includegraphics[width=0.24\linewidth]{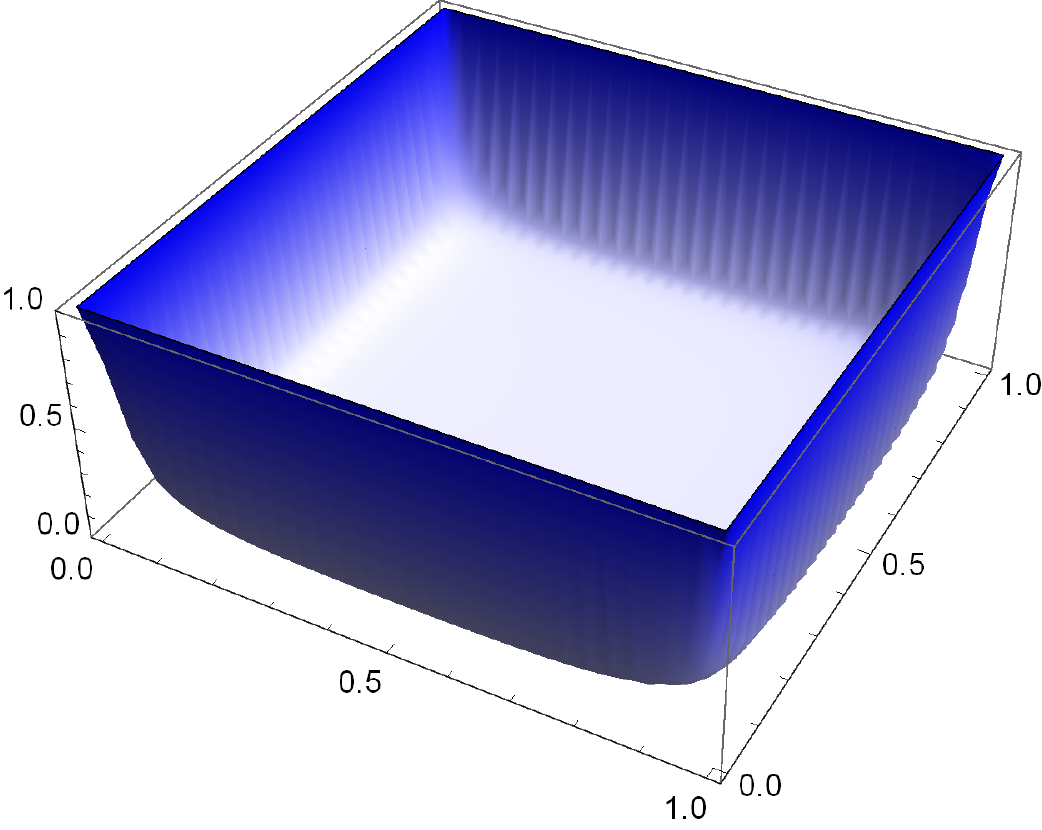}
\includegraphics[width=0.24\linewidth]{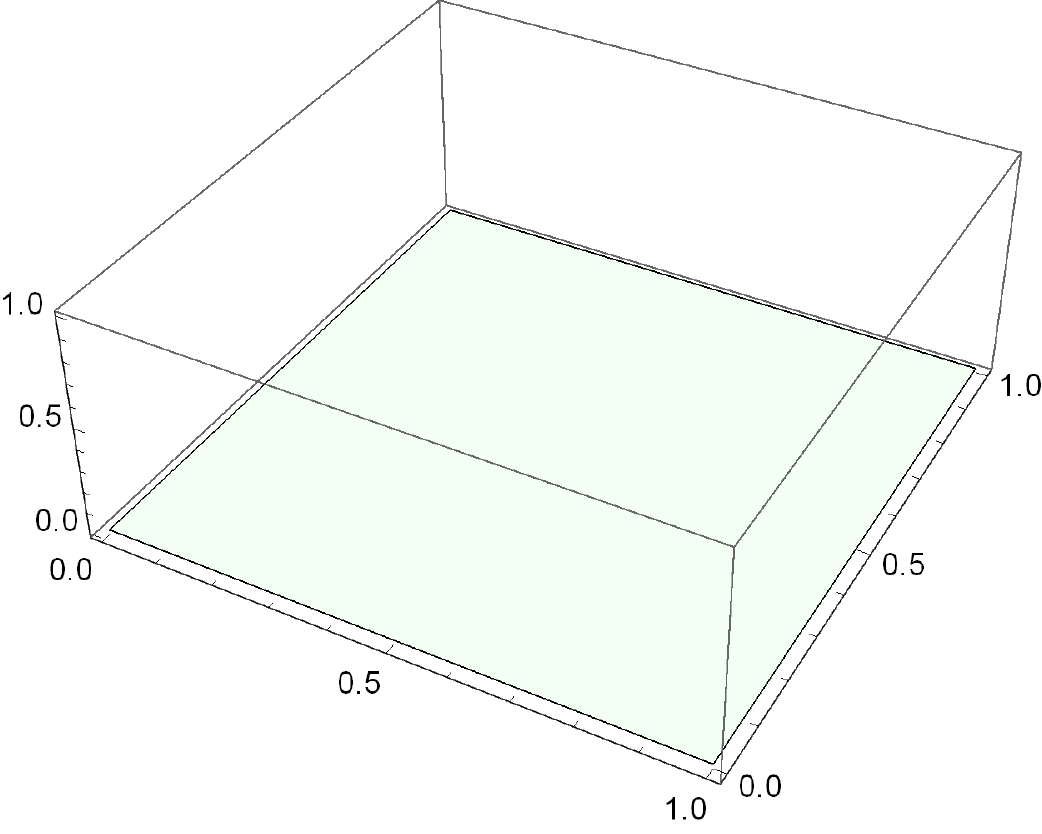}
\includegraphics[width=0.24\linewidth]{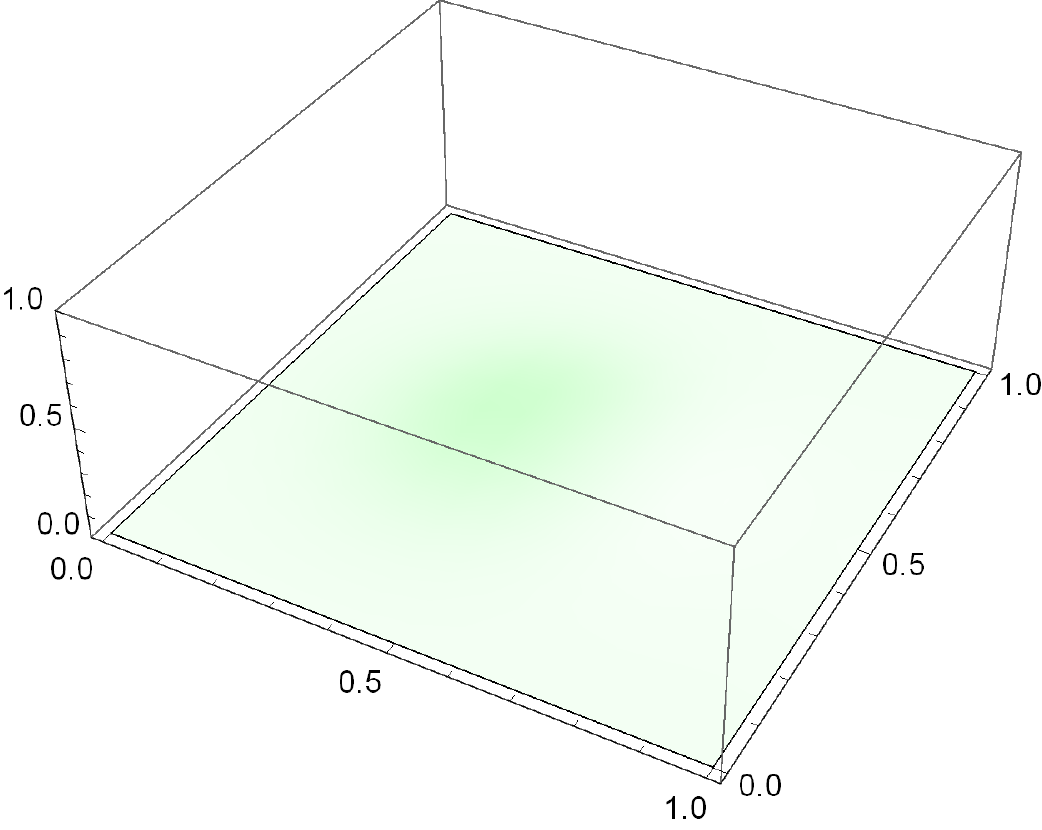}
\includegraphics[width=0.24\linewidth]{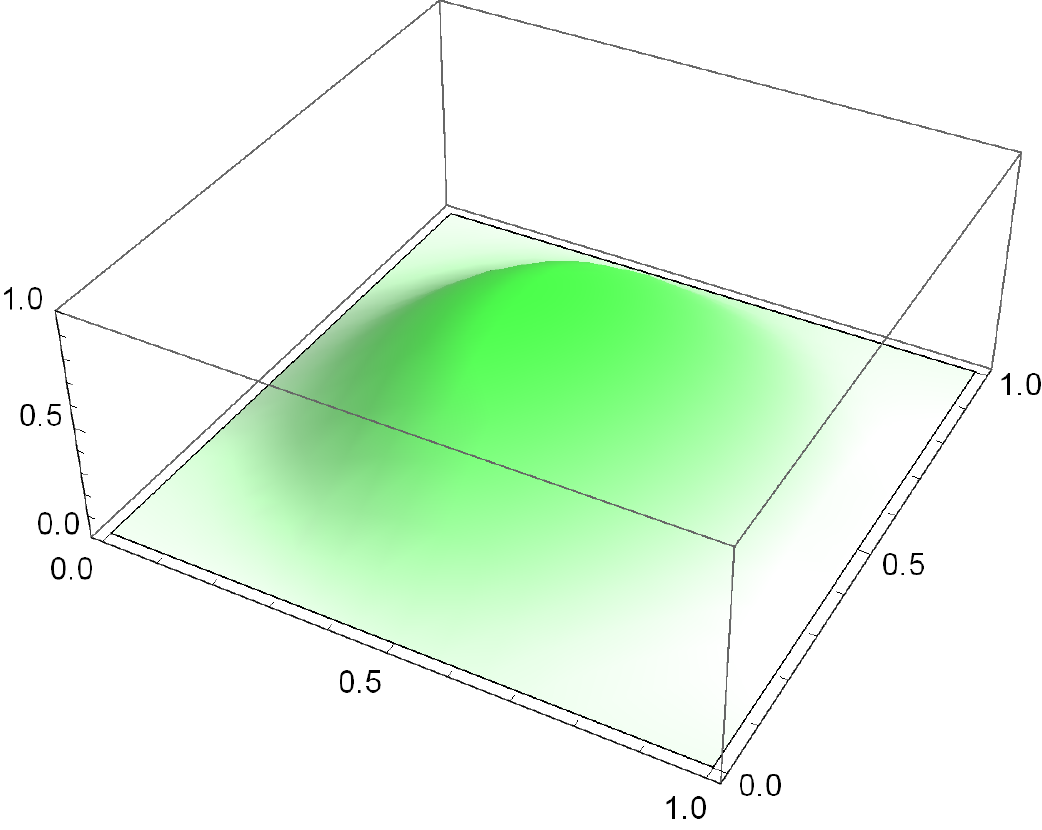}
\includegraphics[width=0.24\linewidth]{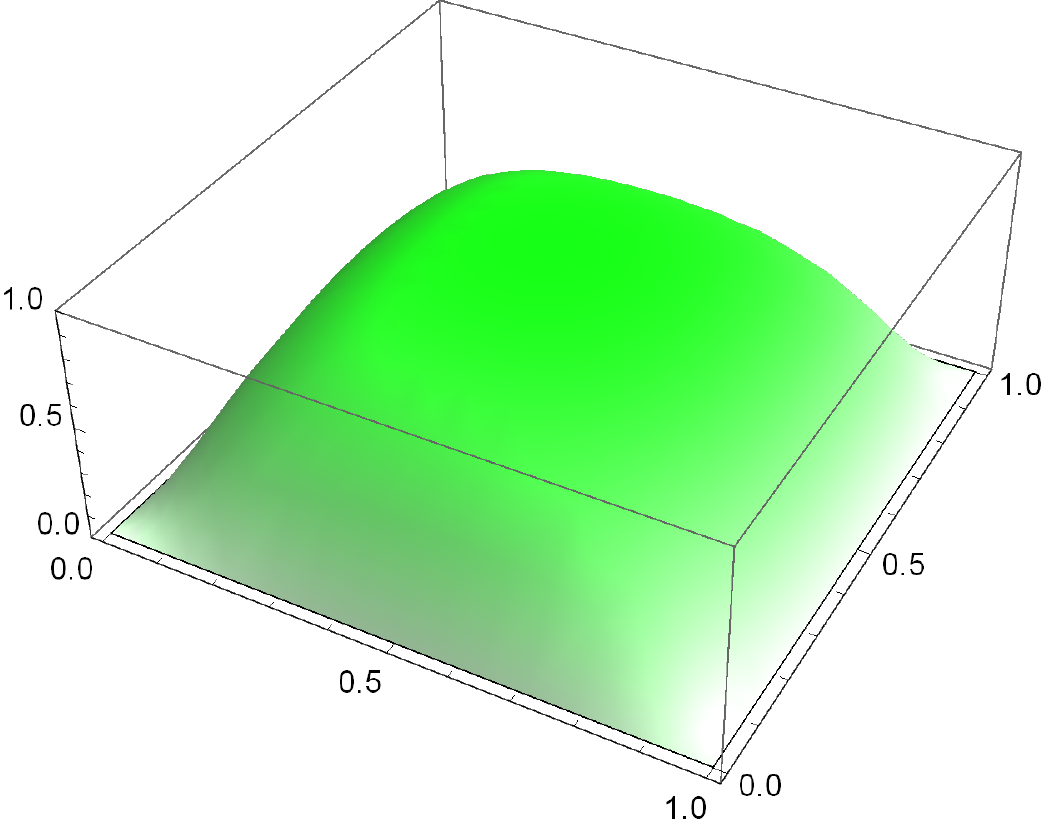}
\caption{Simulation results of model (\ref{0riginalEquations}) with homogeneous Dirichlet boundary conditions. Plots of model solutions $A(x,t)$ (cancer cells, red, top row), $N(x,t)$ (normal cells, blue, middle row) and $H(x,t)$ (lactic acid concentration, green, bottom row) at time points $t=0,33,39,45$, with $x\in \Omega=[0,L]\times [0,L]$. The initial tumor cells survive, produce lactic acid and invade the tissue, leading the drastic reduction on the number of normal cells, except in the boundary, where the normal cells remain at a high number.}
\end{figure}

\newpage

\bibliographystyle{amsplain}

\end{document}